\documentclass[reqno,centertags,12pt,a4paper]{amsart}
\usepackage{dsfont}
\usepackage{amscd}
\usepackage{amsmath,tikz}
\usepackage{amssymb}
\usepackage{amsfonts}
\usepackage{hyperref}
\usepackage{enumerate}
\usepackage{latexsym}
\usepackage{cases}
\usepackage{amsthm}
\usepackage{graphicx}
\usepackage{verbatim}
\usepackage{mathrsfs}
\usepackage{color}

\newcommand\pxx{%
	\mathrel{\text{\tikz[baseline] \draw (0em,-0.3ex) -- (.4em,1.7ex) (.2em,-0.3ex) -- (.6em,1.7ex);}%
}}

\newtheorem{lemma}{Lemma}[section]
\newtheorem{proposition}[lemma]{Proposition}

\newtheorem{definition}[lemma]{Definition}

\newtheorem{remark}[lemma]{Remark}

\newtheorem{theorem}[lemma]{Theorem}
\theoremstyle{definition}
\numberwithin{equation}{section}

\allowdisplaybreaks[2] 

\begin{document}

\title[$L^{p}$-estimate of Schr\"odinger maximal function]{$L^{p}$-estimate of Schr\"odinger maximal function in higher dimensions}
\author[Zhenbin~Cao, ~Changxing ~Miao~and~Meng~Wang]{Zhenbin~Cao, ~Changxing ~Miao~and~Meng~Wang}

\date{\today}

\address{School of Mathematical Sciences, Beijing Normal University, Beijing 100875, China}
\email{11735002@zju.edu.cn}

\address{Institute of Applied Physics and Computational Mathematics, Beijing 100088, China}
\email{miao\_{}changxing@iapcm.ac.cn}

\address{School of Mathematical Science, Zhejiang University, Hangzhou 310058, China}
\email{mathdreamcn@zju.edu.cn}


\begin{abstract}
Almost everywhere convergence on the solution of Schr\"odinger equation is an important problem raised by Carleson in harmonic
analysis. In recent years, this problem was essentially solved by building the sharp $L^p$-estimate of Schr\"odinger maximal function. Du-Guth-Li in \cite{DGL} proved the sharp $L^p$-estimates for all $p \geq 2$ in $\mathbb{R}^{2+1}$. Du-Zhang in \cite{DZ} proved the sharp $L^2$-estimate in $\mathbb{R}^{n+1}$ with $n \geq 3$, but for $p>2$ the sharp $L^p$-estimate of Schr\"odinger maximal function is still unknown. In this paper, we obtain partial results on this problem by using polynomial partitioning.
\end{abstract}

\maketitle

\section{INTRODUCTION}\label{section1}

We consider the free Schr\"{o}dinger equation:
\begin{align}
	\begin{cases}
		\ i u_{t} - \Delta u = 0, \quad (x,t) \in \mathbb{R}^{n} \times \mathbb{R} , \\
		\  u(x,0)=f(x), \ \ \quad  \,  \  x\in \mathbb{R}^{n}.
	\end{cases}
\end{align}
Its solution is given by
\begin{equation*}
	e^{it\Delta}f(x)= (2\pi)^{-n}\int_{\mathbb{R}^{n}} e^{i(x \cdot \xi +t|\xi|^{2})}\widehat{f}(\xi)d\xi ,
\end{equation*}
where $\widehat{f}$ denotes the Fourier transform of the function $f$.

One of the fundamental problems in this setting  is that of determining the  optimal $s$ such that
\begin{equation}\label{aim}
	\lim_{t\rightarrow 0}e^{it\Delta}f(x)=f(x),   \quad \quad \text{for {\rm a.e.}}\, x\in \mathbb R^n
\end{equation}
for all $f\in H^s(\mathbb{R}^n)$. In 1979, Carleson \cite{C} first purposed this problem  and  proved that the  almost everywhere convergence
\eqref{aim} holds for any $f\in H^{1/4}(\mathbb{R})$ by making use of the stationary phase method. A year later, Dahberg and Kenig in \cite{DK}
were able to prove the condition $s\ge\frac14$ given by Carleson is sharp, showing that the existence of functions in $H^s\left(\mathbb{R}\right)$ with $s<\frac14$ for which the
convergence  fails. For the situation in higher dimensions, many authors such as Carbery in \cite{Carbery} and Cowling in \cite{Cowling} studied this problem,
and in 1987 Sj\"{o}lin in \cite{Sjolin} and Vega in \cite{Vega} proved independently that \eqref{aim} holds if
$s>\frac12$ no matter the dimension.
After that some important positive  results in higher dimensions have been obtained by many authors \cite{B1,DG,DGLZ,L,LR,MVV}.
More recently, Bourgain in \cite{B2} gave counterexamples showing that (\ref{aim}) can fail if $s <\frac{n}{2(n+1)}$.
Du-Guth-Li in \cite{DGL} and  Du-Zhang in \cite{DZ} improved the sufficient condition to the almost sharp range $s>\frac{n}{2(n+1)}$
when $n=2$  and $n\ge3$, respectively.  Hence, the Carleson problem  was essentially solved except the endpoint.

Now, we recall the main ideas in \cite{DGL, DZ}. Suppose that $f$ satisfies supp$\widehat{f} \subset B^{n}(0,1)$. When $n=2$, Du-Guth-Li in \cite{DGL} proved the sharp $L^{3}$-estimate of Schr\"odinger maximal function:
\begin{equation}\label{maximal-est-1}
	\Big\|\sup_{0 < t \leq R} |e^{it\Delta}f| \Big\|_{L^{3}(B^{2}(0,R))} \lesssim R^{\epsilon}\|f\|_{L^{2}(\mathbb R^2)}, \;\;\;\forall\, \varepsilon>0.\end{equation}
By Littlewood-Paley decomposition  and parabolic rescaling, they applied (\ref{maximal-est-1}) to derive that (\ref{aim}) holds for $s>1/3$. On the other hand, using H\"{o}lder's inequality, one immediately obtains whole sharp $L^{p}$-estimates:
\begin{equation}\label{maximal-est-2}
	\Big\|\sup_{0 < t \leq R} |e^{it\Delta}f| \Big\|_{L^{p}(B^{2}(0,R))} \lesssim
	\left\{ \begin{aligned}
		& R^{\frac{2}{p}-\frac{2}{3}+\epsilon}\|f\|_{L^{2}(\mathbb R^2)}, & 2 \leq p < 3,\\
		& R^{\epsilon}\|f\|_{L^{2}(\mathbb R^2)}, & p\geq 3.
	\end{aligned}\right.
\end{equation}
This fact means that  $p=3$ is the minimal number such that
\begin{equation}\label{maximal-est-3}\Big\|\sup_{0 < t \leq R} |e^{it\Delta}f| \Big\|_{L^{p}(B^{2}(0,R))} \lesssim R^{\epsilon}\|f\|_{L^{2}(\mathbb R^2)},\;\; \;\;\forall\,\varepsilon>0.\end{equation}
When $n \geq 3$, Du-Zhang in \cite{DZ}  showed (\ref{aim}) holds for $s>\frac{n}{2(n+1)}$  by establishing the sharp $L^{2}$-estimate
\begin{equation}\label{maximal-est-5}\Big\|\sup_{0 < t \leq R} |e^{it\Delta}f| \Big\|_{L^{2}(B^{n}(0,R))} \lesssim R^{\frac{n}{2(n+1)}+\epsilon}\|f\|_{L^{2}(\mathbb R^n)}, \;\;\forall\,\varepsilon>0.\end{equation}
We observe that (\ref{maximal-est-5}) can recover the case of $p=2$ in (\ref{maximal-est-2}) when $n=2$. The sharp $L^p$-estimate of Schr\"odinger maximal function implies the sharp $L^2$-estimate by H\"older's inequality, but the reverse is invalid. So Du-Zhang in \cite{DZ} further raised the problem of identifying the sharp exponent  $p$ such that
\begin{equation}\label{main aim0}
	\left\|\sup_{0 < t \leq R} |e^{it\Delta}f| \right\|_{L^{p}(B^{n}(0,R))} \lesssim R^{\epsilon}\|f\|_{L^{2}(\mathbb R^n)}, \;\;\forall\,\varepsilon>0.
\end{equation}
From (\ref{maximal-est-1}), we see that
\begin{equation}\label{maximal-est-4}\Big\|\sup_{0 < t \leq R} |e^{it\Delta}f| \Big\|_{L^{\frac{2(n+1)}{n}}(B^{n}(0,R))} \lesssim R^{\epsilon}\|f\|_{L^{2}(\mathbb R^n)},\;\;\forall\,\varepsilon>0\end{equation}
holds when $n=1,2$, see \cite{DGL,KPV}. It would be an amazing result if (\ref{maximal-est-4}) holds for higher dimensions, since $p={2(n+1)}/{n}$ corresponds to  the endpoint exponent of $(n+1)$-dimensional restriction conjecture. Unfortunately, Du-Kim-Wang-Zhang in \cite{DKWZ} proved (\ref{maximal-est-4}) fails when $n \geq 3$. More precisely, by researching Bourgain's counterexample in every intermediate dimension in \cite{B2}, they proved that (\ref{main aim0}) fails if
\begin{equation}\label{DKWZ result}
	p < p_0 :=  \max_{\substack{m\in \mathbb{Z}   \\ 1\leq m\leq n}} 2+\frac{4}{n-1+m+\frac{n}{m}}.
\end{equation}
In particular, (\ref{maximal-est-4}) fails for $n \geq 3$ due to $\frac{2(n+1)}{n} <p_0$. Recently,  Wu in \cite{W2} proved that (\ref{main aim0}) holds for
\begin{equation}\label{Wu result}
	p \geq 2+\frac{4}{n+2-\frac{1}{n}}.
\end{equation}

Our main result in this paper is the following:

\begin{theorem}\label{th1}
	Let $n \geq 3$. Suppose that $p$ satisfies
	\begin{equation}\label{main range}
		p\geq 2+\frac{4}{n+1+\frac{1}{2}+...+\frac{1}{n}}.
	\end{equation}
For any $0<\epsilon \ll 1$, there exists a constant $C_{\epsilon}$ such that
\begin{equation}\label{main aim}
	\Big\|  \sup_{0<t\leq R} |e^{it\Delta}f|  \Big\|_{L^{p}(B^{n}(0,R))} \leq C_{\epsilon} R^{\epsilon}\|f\|_{L^{2}}
\end{equation}
for all $f$ with ${\rm supp}\widehat{f} \subset B^{n}(0,1)$.
\end{theorem}



We firstly recall the $k$-broad norm $BL^p_{k,A}$ initially defined by Guth in \cite{G2}. $BL^p_{k,A}$-norm has two advantages compared with classical $L^p$-norm. The first is that $BL^p_{k,A}$-norm can lead to the $k$-linear reduction. Secondly, $\|e^{it\Delta}f\|_{BL^p_{k,A}(B_R)}$ is negligible whenever most of the mass of $e^{it\Delta}f$ is concentrated on an algebraic surface with the dimension less than $k$, named the vanishing property. Guth applied these properties to improve the bound of restriction conjecture in higher dimensions. Motivated by Guth's idea, Du-Li in \cite{DL} extended $BL^p_{k,A}$-norm to a mixed version $BL^{p,\infty}_{k,A}$-norm for studying $L^p$-estimate of Schr\"odinger maximal function. The definition we present below is slightly different from that of their article. They introduced a new parameter $M$ in the definition which can help them to achieve the bilinear reduction in the main proof \cite[Section 2]{DL}. However this makes that  the vanishing property may invalidate when $M$ is large. And it will lead to an obstacle in our following argument. For reserving this important property, we give the following definition. On the bilinear reduction, we are going to take a different way to achieve it.

Let  $f$ satisfy supp$\widehat{f} \subset B^{n}(0,1)$. Decomposing $B^{n}(0,1)$ into balls $\tau$ of radius $K^{-1}$, where $K$ is a large constant, we have $f=\sum_{\tau}f_{\tau}$ such that each $f_{\tau}$ is Fourier supported on $\tau$. Denote
$$G(\tau)=\big\{ G(\xi) :~\xi \in \tau \big\},\;\;\text{where}\;\;  G(\xi) = \tfrac{(-2\xi,1)}{|(-2\xi,1)|}.  $$
If $V \subset \mathbb{R}^{n} \times \mathbb{R}$ is a subspace, we use $\angle (G(\tau),V)$ to be the smallest angle between any non-zero vectors $v\in V$ and $v'\in G(\tau)$.

Next we  decompose $B^{n}(0,R)$ into balls $B_{K}$ of radius $K$, and decompose $[0,R]$ into intervals $I_{K}$ of length $K$. For a given positive integer $A$, we define
\begin{equation}\label{BLp define 1}
	\mu_{e^{it\Delta}f}(B_{K}\times I_{K}):=\min_{V_{1},...,V_{A}} \left(  \max_{\tau \notin V_{a} \mathrm{~for~ all~} a} \int_{B_{K}\times I_{K}}|e^{it\Delta}f_{\tau}(x)|^{p} dxdt  \right),
\end{equation}
where $V_{1},...,V_{A}$ are $(k-1)$-subspaces of $\mathbb{R}^{n+1}$, and $\tau \notin V_{a}$ means that $\angle(G(\tau),V_{a}) > K^{-1}$. For any subset $ U \subset B_{R}^{\ast}:=B^{n}(0,R) \times [0,R] $, define
\begin{equation}\label{BLp define 2}
	\|e^{it\Delta}f\|_{BL_{k,A}^{p,\infty}(U)}^{p} :=\sum_{B_{K}\subset B(0,R)} \max_{I_{K} \subset [0,R]}\frac{|U \cap (B_{K}\times I_{K})|}{|B_{K}\times I_{K}|}\mu_{e^{it\Delta}f}(B_{K}\times I_{K}).
\end{equation}

We reduce Theorem \ref{th1} to the following Theorem \ref{th3}, and we will explain how Theorem \ref{th3} implies Theorem \ref{th1} in Section \ref{reduction}.
\begin{theorem}\label{th3}
	Let $p$ satisfy (\ref{main range}). For any $0<\epsilon\ll 1$, there exist a positive integer $A$ and a constant $C(\epsilon,K)$ such that
	\begin{equation}
		\left\| e^{it\Delta}f \right\|_{BL^{p,\infty}_{2,A}(B_{R}^{\ast})} \leq C(\epsilon,K) R^{\epsilon}\|f\|_{L^{2}}
	\end{equation}
	for all $R\geq 1$ and all $f$ with {\rm supp}$\widehat{f} \subset B^{n}(0,1)$.
\end{theorem}

Since $BL^{p,\infty}_{k,A}$-norm  is not continuous, it is not convenient to use polynomial partitioning. For this purpose,
we substitute $BL^{p,\infty}_{k,A}$-norm with $BL^{p,q}_{k,A}$-norm, defined as
\begin{equation}\label{12345asd}
	\|e^{it\Delta}f\|_{BL_{k,A}^{p,q}(U)}^{p} :=\sum_{B_{K}\subset B(0,R)} \Big[\sum_{I_{K}\subset [0,R]}\Big(\frac{|U \cap (B_{K}\times I_{K})|}{|B_{K}\times I_{K}|}\mu_{e^{it\Delta}f}(B_{K}\times I_{K})\Big)^{\frac{q}{p}}\Big]^{\frac{p}{q}}.
\end{equation}
Since
$$ \big\| e^{it\Delta}f \big\|_{BL^{p,\infty}_{k,A}(B_{R}^{\ast})} =\lim_{q\rightarrow \infty}\big\| e^{it\Delta}f \big\|_{BL^{p,q}_{k,A}(B_{R}^{\ast})},$$
it suffices to show the following result.

\begin{theorem}\label{th4}
	Let $p$ satisfy (\ref{main range}). For any $0<\epsilon\ll 1$, let $0<\delta\ll \epsilon$, there exist a positive integer
	$A$ and a constant $C(\epsilon,K)$ such that for any $q>1/\delta$,
	\begin{equation}
		\| e^{it\Delta}f \|_{BL^{p,q}_{2,A}(B_{R}^{\ast})} \leq C(\epsilon,K) R^{\epsilon}\|f\|_{L^{2}}
	\end{equation}
	for all $R\geq 1$ and all $f$ with {\rm supp}$\widehat{f} \subset B^{n}(0,1)$.
\end{theorem}

Let us say a few words on the proof  scheme of Theorem \ref{th4}. The basic tool is polynomial partitioning introduced by Guth in \cite{G2}, which helps us to make the dimensional reduction. We will adopt the algorithm, repeated application of polynomial partitioning, which has similar spirit as Hickman-Rogers in \cite{HR1}.
Hickman-Rogers built the algorithm to estimate $\|Ef\|_{BL^{p}_{k,A}(B_R)}$, where $Ef$ represents the extension operator associated with restriction problem. By building polynomial partitioning on $BL^{p,q}_{k,A}$-norm, we show the algorithm also works to estimate $\|e^{it\Delta}f\|_{BL^{p,q}_{k,A}(B_R^\ast)}$.

This paper is organized as follows. In Section \ref{reduction}, we explain how Theorem \ref{th3} implies Theorem \ref{th1}. In Section \ref{section3}, we introduce some notations and standard results. The proof of Theorem \ref{th4} is in Section \ref{section4} and Section \ref{section5}.
\vskip0.2cm

To end up this section, we will outline some notations used throughout the paper.
If $X$ is a finite set, we use $|X|$ to denote its cardinality. If $X$ is a measurable set,
we use $|X|$ to denote its Lebesgue measure. If the function $f$ has compact support,
we use ${\rm supp}f$ to denote the support of $f$. $C_{\epsilon}$  denotes a constant which depends on
$\epsilon$.  Write $A\lesssim B$ or $A=O(B)$ to
mean that there exists a constant $C$ such that $A\leq CB$. We write $ A \lessapprox B$ to denote $A \lesssim_{\epsilon} R^{\epsilon} B$.

\section{How maximal broad $L^{p}$-norm implies maximal $L^{p}$-norm}\label{reduction}

In this section, we prove that Theorem \ref{th3} implies Theorem \ref{th1} by broad-narrow argument. Broad-narrow argument was developed by Bourgain and Guth \cite{BG11} to study restriction conjecture for the first time. This method can help us to reduce original linear estimate to $k$-linear estimate. Though the broad norm that we consider in the paper is slightly weaker than classical $k$-linear norm, the method and reduction still work.

For the convenience of the discussion, we first introduce an equivalent definition of $BL^{p,q}_{k,A}$-norm. Denote
$$  Br_{k,A}^K(e^{it\Delta}f)(x):=\min_{V_{1},...,V_{A}} \max_{\tau \notin V_{a} \mathrm{~for~ all~} a}  |e^{it\Delta}f_\tau(x)|,   $$
where $V_1$,...,$V_A$ are defined as in Section 1. Then we have
\begin{equation}\label{new f}
\|e^{it\Delta}f\|_{BL_{k,A}^{p,q}(B_R^\ast)}\sim \Big\|Br_{k,A}^K(e^{it\Delta}f)(x)\Big\|_{L_x^p L_t^q(B_R^\ast)}.
\end{equation}
In fact, note that $\tau$ is a $K^{-1}$-ball, we can view $|e^{it\Delta}f_{\tau}|$ as constant on $B_K \times I_K$. So (\ref{BLp define 1}) becomes
\begin{equation*}
	\mu_{e^{it\Delta}f}(B_{K}\times I_{K})\sim  \int_{B_{K}\times I_{K}} \min_{V_{1},...,V_{A}}  \max_{\tau \notin V_{a} \mathrm{~for~ all~} a}     |e^{it\Delta}f_{\tau}(x)|^{p} dxdt .
\end{equation*}
Summing all $B_K\subset B(0,R)$ and all $I_K \subset [0,R]$ in (\ref{12345asd}), we can get (\ref{new f}) immediately. Here we have used one important property named locally constant property. We will introduce it more precisely in Section 3.

\begin{proposition}
For every $0<\epsilon \ll1 $ and every sufficiently large $R$, there exist
$$   1\ll A \ll K_1 \ll ...\ll K_{n-1} \ll K_n \ll R^\epsilon   $$
such that the following holds.
Suppose that $f$ is Fourier supported on $B^n(0,1)$, then
	\begin{align}
		|e^{it\Delta}f(x)|&\lesssim K_n^{2n} \sum_{\alpha_1,...,\alpha_{n+1} \text{transverse}} \Big( \prod_{j=1}^{n+1} |e^{it\Delta}f_{\alpha_j}(x)| \Big)^{\frac{1}{n+1}}+K_{n-1}^{2(n-1)} \max_{V^n} Br_{n,A}^{K_{n-1}}(e^{it\Delta}f_{V^n})(x) \nonumber\\
		&+ A K_{n-2}^{2(n-2)} \max_{V^{n-1}} Br_{n-1,A}^{K_{n-2}}(e^{it\Delta}f_{V^{n-1}})(x)+...+A^{n-2} K_{1}^{2} \max_{V^{2}} Br_{2,A}^{K_{1}}(e^{it\Delta}f_{V^{2}})(x)  \nonumber\\
		&+A^{n-1} \max_{\gamma:K_1^{-1}-\text{balls}}|e^{it\Delta}f_\gamma(x)|.\label{bil red1}
	\end{align}
Here each $\alpha_i$ denotes the ball with radius $K_n^{-1}$, and each $V^l$ denotes one $l$-dimensional subspace. The definition of transverse see (\ref{trans cond se2}).
And for each $1\leq l \leq n$, we have
	\begin{align}
		|e^{it\Delta}f(x)|&\lesssim K_{l-1}^{2(l-1)} \max_{V^l} Br_{l,A}^{K_{l-1}}(e^{it\Delta}f_{V^l})(x) + A K_{l-2}^{2(l-2)} \max_{V^{l-1}} Br_{l-1,A}^{K_{l-2}}(e^{it\Delta}f_{V^{l-1}})(x)+... \nonumber\\
		&+A^{n-2} K_{1}^{2} \max_{V^{2}} Br_{2,A}^{K_{1}}(e^{it\Delta}f_{V^{2}})(x)  +A^{n-1} \max_{\gamma:K_1^{-1}-\text{balls}}|e^{it\Delta}f_\gamma(x)|.\label{bil red2}
	\end{align}
	
\end{proposition}

\noindent\textit{Proof.}  We only prove (\ref{bil red1}) since the proof of (\ref{bil red2}) is the same. We first fix $x$ and $t$.
Decomposing the frequency support $B^n(0,1)$ to balls $\alpha$ with radius $K_n^{-1}$, then $$e^{it\Delta}f(x)=\sum_\alpha e^{it\Delta}f_\alpha(x),$$
where each $f_\alpha$ is Fourier supported in $\alpha$. Define
$$  \mathcal{S}(x,t):=\Big\{ \alpha :|e^{it\Delta}f_\alpha(x)|\geq \frac{1}{K_n^n} \max_{\alpha'} |e^{it\Delta}f_{\alpha'}(x)| \Big\}.  $$
There are two cases that may occur now:
\begin{itemize}
	\item There exist $\alpha_1,...,\alpha_{n+1}\in \mathcal{S}(x,t) $ which are $(n+1)$-transverse: for any $v_j \in G(\alpha_j)$,
	\begin{equation}\label{trans cond se2}
		|v_1 \wedge v_2\wedge...\wedge v_{n+1}| \gtrsim K_{n}^{-n}.
	\end{equation}
	\item There exists an $n$-dimensional subspace $V^n$ such that for each $\alpha \in  \mathcal{S}(x,t),$
	$$\angle(G(\alpha),V^n)\leq \frac{1}{K_n}.$$
\end{itemize}
If we are in the first case,
$$   |e^{it\Delta}f(x)|\leq K^n_n \max_{\alpha'} |e^{it\Delta}f_{\alpha'}(x)|\leq K_n^{2n} \Big( \prod_{j=1}^{n+1}|e^{it\Delta}f_{\alpha_j}(x)| \Big)^{\frac{1}{n+1}}.  $$
If we are in the second case,
\begin{align*}
	|e^{it\Delta}f(x)|&=\Big|\sum_{\alpha\in  \mathcal{S}(x,t)} e^{it\Delta}f_\alpha(x) + \sum_{\alpha\notin  \mathcal{S}(x,t)}e^{it\Delta}f_\alpha(x) \Big|   \\
	& \leq \Big|\sum_{\alpha\in  \mathcal{S}(x,t)} e^{it\Delta}f_\alpha(x) \Big|+\max_{\alpha'} |e^{it\Delta}f_{\alpha'}(x)|  \\
	& \lesssim \Big|  \sum_{\substack{ \alpha: \\ \angle(G(\alpha),V^n)\leq K_n^{-1}  }} e^{it\Delta}f_\alpha(x) \Big| +\max_{\alpha'}|e^{it\Delta}f_{\alpha'}(x)| .
\end{align*}
Therefore we conclude
$$  |e^{it\Delta}f(x)|\lesssim K_n^{2n} \Big( \prod_{j=1}^{n+1}|e^{it\Delta}f_{\alpha_j}(x)| \Big)^{\frac{1}{n+1}}+  \Big|  \sum_{\substack{ \alpha: \\ \angle(G(\alpha),V^n)\leq K_n^{-1}  }} e^{it\Delta}f_\alpha(x) \Big| +\max_{\alpha'}|e^{it\Delta}f_{\alpha'}(x)|.   $$

Denote
$$ e^{it\Delta}f_{V^n}(x):= \sum_{\substack{ \alpha: \\ \angle(G(\alpha),V^n)\leq K_n^{-1}  }} e^{it\Delta}f_\alpha(x).  $$
Decomposing the frequency support $N_{K^{-1}_n}(V^n)$ to balls $\beta$ with radius $K_{n-1}^{-1}$, then
$$e^{it\Delta}f_{V^n}(x)=\sum_{\substack{\beta:   \\  G(\beta) \cap N_{K^{-1}_n}(V^n)\neq \varnothing }} e^{it\Delta}f_{V^n,\beta}(x).$$
Define
$$  \mathcal{S}_{V^n}(x,t):=\Big\{ \beta :|e^{it\Delta}f_{V^n,\beta}(x)|\geq \frac{1}{K_{n-1}^{n-1}} \max_{\beta'} |e^{it\Delta}f_{V^n,\beta'}(x)|,   G(\beta) \cap N_{K^{-1}_n}(V^n)\neq \varnothing\Big\}.  $$
We choose $(n-1)$-dimensional subspaces $V^{n-1}_1,...,V^{n-1}_A\subset V^n$ to achieve the minimum of
$$   \max_{ \substack{\beta:  \\  \angle(G(\beta),V_a^{n-1})>K^{-1}_{n-1}, \forall a   } }|e^{it\Delta}f_{V^n,\beta}(x)| . $$
There are two cases  that may occur now:
\begin{itemize}
	\item There exists $\beta \in \mathcal{S}_{V^n}(x,t) $ such that for each $a$, $\angle( G(\beta),V_a^{n-1}  )>K^{-1}_{n-1}.$
	\item For each $\beta \in \mathcal{S}_{V^n}(x,t) $, there exists  $a$ such that $\angle( G(\beta),V_a^{n-1}  )\leq K^{-1}_{n-1}.$
\end{itemize}
If the first case happens,
\begin{align*}
	|e^{it\Delta}f_{V^n}(x)| &\leq K^{n-1}_{n-1} \max_{\beta'}|e^{it\Delta}f_{V^n,\beta'}(x)| \\
	& \leq K^{2(n-1)}_{n-1} \min_{\beta \in \mathcal{S}_{V^n}(x,t)} |e^{it\Delta}f_{V^n,\beta}(x)| \\
	&  \leq K^{2(n-1)}_{n-1} \max_{ \substack{\beta:  \\  \angle(G(\beta),V_a^{n-1})>K^{-1}_{n-1}, \forall a   }} |e^{it\Delta}f_{V^n,\beta}(x)|    \\
	&=K^{2(n-1)}_{n-1} Br_{n,A}^{K_{n-1}}(e^{it\Delta}f_{V^n})(x).
\end{align*}
If the second case happens,
\begin{align*}
	|e^{it\Delta}f_{V^n}(x)|&=\Big|\sum_{\beta \in \mathcal{S}_{V^n}(x,t)} e^{it\Delta}f_{V^n,\beta}(x) + \sum_{\beta\notin  \mathcal{S}_{V^n}(x,t)}e^{it\Delta}f_{V^n,\beta}(x) \Big|   \\
	& \leq  \Big|\sum_{\beta \in \mathcal{S}_{V^n}(x,t)} e^{it\Delta}f_{V^n,\beta}(x) \Big|+\max_{\beta'} |e^{it\Delta}f_{V^n,\beta'}(x)|  \\
	& \leq \sum_{a=1}^A \Big|  \sum_{\substack{\beta \in \mathcal{S}_{V^n}(x,t) \\ \angle(G(\beta),V_a^{n-1})\leq K_{n-1}^{-1}  }} e^{it\Delta}f_{V^n,\beta}(x) \Big| +\max_{\beta'} |e^{it\Delta}f_{V^n,\beta'}(x)|    \\
	& \leq \sum_{a=1}^A \Big|  \sum_{\substack{\beta \in \mathcal{S}_{V^n}(x,t) \\ \angle(G(\beta),V_a^{n-1})\leq K_{n-1}^{-1}  }} e^{it\Delta}f_{\beta}(x) \Big|    \\
	&\quad + \sum_{a=1}^A \Big|  \sum_{\substack{\beta \in \mathcal{S}_{V^n}(x,t) \\ \angle(G(\beta),V_a^{n-1})\leq K_{n-1}^{-1}  }} \Big(e^{it\Delta}f_{V^n,\beta}(x)-e^{it\Delta}f_\beta(x)\Big) \Big|  +\max_{\beta'} |e^{it\Delta}f_{V^n,\beta'}(x)|   \\
	&= \sum_{a=1}^A \Big|  \sum_{\substack{\beta \in \mathcal{S}_{V^n}(x,t) \\ \angle(G(\beta),V_a^{n-1})\leq K_{n-1}^{-1}  }} e^{it\Delta}f_{\beta}(x) \Big|    \\
	&\quad + \sum_{a=1}^A \Big|  \sum_{\substack{\beta \in \mathcal{S}_{V^n}(x,t) \\ \angle(G(\beta),V_a^{n-1})\leq K_{n-1}^{-1}  }}\sum_{  \substack{\alpha \subset \beta \\ \angle( G(\alpha), V^n )> K^{-1}_n}  } e^{it\Delta}f_\alpha(x) \Big|  +\max_{\beta'} |e^{it\Delta}f_{V^n,\beta'}(x)|   \\
	& \leq \sum_{a=1}^A \Big|  \sum_{\substack{\beta \in \mathcal{S}_{V^n}(x,t) \\ \angle(G(\beta),V_a^{n-1})\leq K_{n-1}^{-1}  }} e^{it\Delta}f_{\beta}(x) \Big| +  A\sum_{\alpha \notin \mathcal{S}(x,t)} |e^{it\Delta}f_\alpha(x)|  + \max_{\beta'} |e^{it\Delta}f_{V^n,\beta'}(x)|    \\
	& \lesssim   \sum_{a=1}^A \Big|  \sum_{\substack{\beta\\ \angle(G(\beta),V_a^{n-1})\leq K_{n-1}^{-1}  }} e^{it\Delta}f_{\beta}(x) \Big| +  A\max_{\alpha'} |e^{it\Delta}f_{\alpha'}(x)|  + \max_{\beta'} |e^{it\Delta}f_{V^n,\beta'}(x)|.
\end{align*}
Thus
$$   |e^{it\Delta}f_{V^n}(x)|  \lesssim   K^{2(n-1)}_{n-1} Br_{n,A}^{K_{n-1}}(e^{it\Delta}f_{V^n})(x)  +    \sum_{a=1}^A \Big|  \sum_{\substack{\beta\\ \angle(G(\beta),V_a^{n-1})\leq K_{n-1}^{-1}  }} e^{it\Delta}f_{\beta}(x) \Big| +  A\max_{\alpha'} |e^{it\Delta}f_{\alpha'}(x)|  .  $$

Combining the results of the first two steps, we obtain
\begin{align*}
	|e^{it\Delta}f(x)| &  \lesssim K_n^{2n} \Big( \prod_{j=1}^{n+1}|e^{it\Delta}f_{\alpha_j}(x)| \Big)^{\frac{1}{n+1}}  \\
	&\quad +   K^{2(n-1)}_{n-1} Br_{n,A}^{K_{n-1}}(e^{it\Delta}f_{V^n})(x)  +    \sum_{a=1}^A \Big|  \sum_{\substack{\beta\\ \angle(G(\beta),V_a^{n-1})\leq K_{n-1}^{-1}  }} e^{it\Delta}f_{\beta}(x) \Big|   \\
	& \quad  + A\max_{\alpha'} |e^{it\Delta}f_{\alpha'}(x)|  .
\end{align*}
Iterating on this formula, we can get
\begin{align*}
	|e^{it\Delta}f(x)|&\lesssim K_n^{2n}  \Big( \prod_{j=1}^{n+1} |e^{it\Delta}f_{\alpha_j}(x)| \Big)^{\frac{1}{n+1}}+K_{n-1}^{2(n-1)}  Br_{n,A}^{K_{n-1}}(e^{it\Delta}f_{V^n})(x) \nonumber\\
	&+ A K_{n-2}^{2(n-2)}  Br_{n-1,A}^{K_{n-2}}(e^{it\Delta}f_{V^{n-1}})(x)+...+A^{n-2} K_{1}^{2}  Br_{2,A}^{K_{1}}(e^{it\Delta}f_{V^{2}})(x)  \nonumber\\
	&+A^{n-1} \max_{\gamma:K_1^{-1}-\text{balls}}|e^{it\Delta}f_\gamma(x)|.
\end{align*}
Finally, since the choices of $\alpha_j$, $V^2$,..., $V^n$ depend on $x$ and $t$, we take the maximum at both ends of above inequality simultaneously, and then (\ref{bil red1}) follows.

\qed

Now we show that Theorem \ref{th3} implies Theorem \ref{th1}. Take $l=2$ in (\ref{bil red2}), then
$$		|e^{it\Delta}f(x)|\lesssim A^{n-2} K_{1}^{2} \max_{V^{2}} Br_{2,A}^{K_{1}}(e^{it\Delta}f_{V^{2}})(x)  +A^{n-1} \max_{\gamma:K_1^{-1}-\text{balls}}|e^{it\Delta}f_\gamma(x)|.   $$
We integrate over $B_R^\ast$ to obtain
\begin{align*}
	\Big\|  \sup_{0<t\leq R} |e^{it\Delta}f|  \Big\|_{L^{p}(B^{n}(0,R))} &\lesssim  A^{n-2} K_{1}^{2} \Big\|  \max_{V^{2}} Br_{2,A}^{K_{1}}(e^{it\Delta}f_{V^{2}})(x) \Big\|_{L_x^p L_t^\infty(B_R^\ast)}  \\
	& +   A^{n-1} \sum_{\gamma}    \Big(	\Big\|  \sup_{0<t\leq R} |e^{it\Delta}f_\gamma|  \Big\|^p_{L^{p}(B^{n}(0,R))}  \Big)^{\frac{1}{p}} .
\end{align*}
For the first term, we have
$$ A^{n-2} K_{1}^{2} \Big\|  \max_{V^{2}} Br_{2,A}^{K_{1}}(e^{it\Delta}f_{V^{2}})(x) \Big\|_{L_x^p L_t^\infty(B_R^\ast)} \leq  A^{n-2} K_{1}^{2}   \sum_{V^{2}}\Big\|  Br_{2,A}^{K_{1}}(e^{it\Delta}f_{V^{2}})(x) \Big\|_{L_x^p L_t^\infty(B_R^\ast)}. $$
Note that the choice of $V^2$ is allowed a $K_2^{-1}$-scale of perturbation, so we can only consider $K_2^{n}$ many $V^2$. Then using (\ref{new f}) and Theorem \ref{th3}, it is further bounded by
$$  \sim A^{n-2} K_{1}^{2}   \sum_{V^{2}} \|e^{it\Delta}f_{V^2}\|_{BL_{2,A}^{p,\infty}(B_R^\ast)} \lesssim    A^{n-2} K_{1}^{O(1)}  K_2^n \|f\|_{L^2}  \lesssim R^\epsilon \|f\|_{L^2}.   $$
For the second term, by parabolic rescaling and induction on scale, one gets
$$  A^{n-1} \sum_{\gamma}    \Big(	\Big\|  \sup_{0<t\leq R} |e^{it\Delta}f_\gamma|  \Big\|^p_{L^{p}(B^{n}(0,R))}  \Big)^{\frac{1}{p}} \leq  A^{n-1} K_1^{-\frac{n}{2}+\frac{n}{p}-\epsilon} R^\epsilon \Big(\sum_\gamma \|f_\gamma\|_{L^2}^p\Big)^{\frac{1}{p}}. $$
Since $p \geq 2$, we apply Minkowski's inequality to obtain
$$ A^{n-1} \sum_{\gamma}    \Big(	\Big\|  \sup_{0<t\leq R} |e^{it\Delta}f_\gamma|  \Big\|^p_{L^{p}(B^{n}(0,R))}  \Big)^{\frac{1}{p}} \leq  A^{n-1} K_1^{-\frac{n}{2}+\frac{n}{p}-\epsilon} R^\epsilon \|f\|_{L^2}.  $$
Finally, we can choose the parameters $A$ and $K_1$ satisfying $A^{n-1} \ll K_1$, and then the induction is close.

Since we will mainly consider $BL^{p,q}_{2,A}$-norm in the following discussion,  from now on,  we use $BL^{p,q}_{A}$ to represent $BL^{p,q}_{2,A}$.

\section{Some preparations}\label{section3}

The main ingredient of the proof of Theorem \ref{th4} is the algorithm which leads to different scales $r$ satisfying $1 \ll r \leq R$, not just the scale $R$. In this section, we introduce some notations and standard results with respect to the general scale $r$.

\vskip0.20cm
\noindent 1. \textbf{Locally constant property}
\vskip0.15cm

Locally constant property says that if a function $f$ has compact Fourier support $\Theta$, then we can view $|f|$ essentially as constant on dual $\Theta^\ast$.

\begin{lemma}\cite[Lemma 6.1]{GWZ}\label{lcp th}
	Let $\Theta$ be a compact symmetric convex set centered at $C_{\Theta} \in \mathbb{R}^{n}$. If $\widehat{g}_{\Theta}$ is supported in $\Theta$ and $T_{\Theta}$ is the dual convex $\Theta^\ast:=\{ y: |y\cdot x- C_{\Theta}|\leq 1, ~\forall ~x\in \Theta  \}$, then there exists a positive function $\eta_{T_{\Theta}}$ satisfying the following properties:
	\begin{itemize}
		\item[{\rm (a)}]\, $\eta_{T_{\Theta}}$ is essentially supported on $10T_{\Theta}$ and rapidly decaying away from it: for any integer $N>0$, there exists a constant $C_N$ such that $\eta_{T_{\Theta}}(x) \leq C_N (1+n(x,10T_\Theta))^{-N}$ where $n(x,10T_\Theta)$ is the smallest positive integer $n$ such that $x\in 10nT_\Theta$,
		\item[{\rm (b)}]\, $\|\eta_{T_{\Theta}}\|_{L^1} \lesssim 1$,
		\item[{\rm (c)}]\, \begin{equation}\label{lcp}
			|g_\Theta|\leq \sum_{T\pxx T_\Theta} c_T \chi_T \leq |g_\Theta|\ast \eta_{T_{\Theta}},
		\end{equation}
		where $c_T:=\max_{x\in T} |g_\Theta|(x)$ and the sum $\sum_{T\pxx T_\Theta}$ is over a finitely overlapping cover $\{T\}$ of $\mathbb{R}^n$ with each $T\pxx T_\Theta$. Here $T\pxx T_\Theta$ means that $T$ is a translated copy of $T_\Theta$.
	\end{itemize}
\end{lemma}

In particular, for each ball $\tau$ of radius $K^{-1}$, note $e^{it\Delta}f_\tau$ is Fourier supported on $\{ (w,|w|^2):w\in \tau  \}$  in a distributional sense, which is contained in a rectangular box of radius $K^{-1}$ and thickness $K^{-2}$. Therefore $e^{it\Delta}f_\tau$ is locally constant on each tube of radius $K$ and length $K^{2}$.

\vskip0.20cm
\noindent 2. \textbf{Wave packet decomposition on scale $r$}
\vskip0.15cm

Let $\varphi$ be a Schwartz function satisfying supp$\widehat{\varphi} \subset B(0,3/2)$, and
$$\sum_{k\in \mathbb{Z}^{n}} \widehat{\varphi}(\xi -k)=1,\quad  \text{~for ~all~}  \xi\in \mathbb{R}^{n}.$$
Define
\begin{equation*}\left\{\begin{aligned}&\widehat{\varphi_{\theta}}(\xi):=r^{\frac{n}{2}}\widehat{\varphi}\left(r^{\frac{1}{2}}(\xi-c(\theta))\right),~ \widehat{\varphi_{\theta,\nu}}(\xi):=e^{-ic(\nu)\cdot \xi}\widehat{\varphi_{\theta}}(\xi),\\
		&c(\theta) \in r^{-1/2}\mathbb{Z}^{n}, \quad\;\; c(\nu) \in r^{1/2}\mathbb{Z}^{n},\end{aligned}\right.\end{equation*}
where $\theta$ denote $r^{-\frac12}$-balls in frequency space  and $\nu$ denote  $r^{\frac12}$-balls in  physical space.

For each Schwartz function $f$ with ${\rm supp}\widehat{f}\subset B^n(0,1)$, we have the following decomposition
\begin{equation}\label{wp 1}
	f  = \sum_{\theta,\nu} c_{\theta,\nu}\varphi_{\theta,\nu},\qquad\; \theta \cap B^{n}(0,1)\not=\varnothing,
\end{equation}
where each $\varphi_{\theta,\nu}$ is Fourier supported in $\theta$ and has physical support essentially in $\nu$. A basic property is the functions $\varphi_{\theta,\nu}$ are approximately orthogonal, i.e.
\begin{equation}\label{wp 2}
	\sum_{\theta,\nu}|c_{\theta,\nu}|^{2} \sim \|f\|_{L^{2}}^{2}.
\end{equation}

Next we consider  the decomposition on $e^{it\Delta}f$ associated with \eqref{wp 1}
\begin{equation}\label{wp 3}
	e^{it\Delta}f = \sum_{\theta,\nu}c_{\theta,\nu}e^{it\Delta}\varphi_{\theta,\nu}.
\end{equation}
By the stationary phase method, we obtain
\begin{equation}\label{wp 4}
	\left| e^{it\Delta}\varphi_{\theta,\nu}(x)   \right| \leq r^{-\frac{n}{4}}\chi_{T_{\theta,\nu}}(x,t) + O(r^{-N}), \;\;\forall \; N>0.
\end{equation}
Here $T_{\theta,\nu}$ is defined by
$$ T_{\theta,\nu} := \left\{(x,t) \in \mathbb{R}^{n} \times \mathbb{R} :0 \leq t \leq r,~ \left|x + 2tc(\theta)-c(\nu) \right| \leq r^{\frac{1}{2}+\delta}  \right\},$$
where $\delta$ is a small positive parameter satisfying $\delta \ll \epsilon.$ $T_{\theta,\nu}$ is a tube of length $r$, of radius $r^{1/2+\delta}$ and in the direction $(-2c(\theta),1)$. For each $\nu$, $e^{it\Delta}\varphi_{\theta,\nu}$ is Fourier supported in $\xi(\theta)=\{ (\xi,\xi^{2}):\xi \in \theta \}$ in a distributional sense. For detailed proofs of wave packet decomposition  one can see \cite{DGL,T}.

\vskip0.25cm

\noindent 3. \textbf{The properties of} $BL^{p,q}_{A}$\textbf{-norm}

\vskip0.20cm
Firstly, we recall  the basic properties of $BL^{p,q}_{A}$-norm.

\begin{proposition}[\cite{DL,G2}]\label{BLp pro} Let $p \leq q$. The following properties hold:
	\begin{itemize}
		\item[{\rm (a)}]\,{\rm (Finite subadditivity)} Let $A\ge1$,  $U_{1},U_{2}\subset \mathbb{R}^{n+1}$, then
		$$   \| e^{it\Delta}f \|^{p}_{BL^{p,q}_{A}(U_{1}\cup U_{2})} \leq \| e^{it\Delta}f \|^{p}_{BL^{p,q}_{A}(U_{1})}+\| e^{it\Delta}f \|^{p}_{BL^{p,q}_{A}(U_{2})}.  $$
		\item[{\rm (b)}]\,{\rm (Triangle inequality)} Let $A \geq 2$, $U\subset \mathbb{R}^{n+1}$, then
		$$   \| e^{it\Delta}(f+g) \|^{p}_{BL^{p,q}_{A}(U)} \lesssim \| e^{it\Delta}f \|^{p}_{BL^{p,q}_{A/2}(U)}+\| e^{it\Delta}g \|^{p}_{BL^{p,q}_{A/2}(U)}.  $$
		\item[{\rm (c)}]\,{\rm (Logarithmic convexity)} Let $A \geq 2$, $U\subset \mathbb{R}^{n+1}$. Suppose that $1\leq p, p_{1}, p_{2} <\infty$  and $0 \leq \alpha\leq 1$ obey
		$$     \frac{1}{p}=\frac{1-\alpha}{p_{1}} +\frac{\alpha}{p_{2}},  $$
		then
		$$  \| e^{it\Delta}f \|_{BL^{p,q}_{A}(U)} \leq \| e^{it\Delta}f \|_{BL^{p_{1},q}_{A/2}(U)}^{1-\alpha}    \| e^{it\Delta}f \|_{BL^{p_{2},q}_{A/2}(U)}^{\alpha}.       $$
	\end{itemize}
\end{proposition}

The proof of Proposition \ref{BLp pro} is elementary, one can refer to \cite{DL,G2}. If $A=1$, we can't use (b) and (c) in Proposition \ref{BLp pro} since $BL^{p,q}_{A}$-norm only makes sense for the positive integer $A$.
But if we start with a large $A$, and use (b) and (c) only $O_{\epsilon}(1)$ times, then the choice of $A$ won't influence our argument, see \cite{HR1,W1}. On the other hand, we always choose the parameter $A$ satisfying $A \ll K$, which is an unimportant parameter. From now on, we write $BL^{p,q}$ to represent $BL^{p,q}_{A}$.

Now  we consider the relationship between $BL^{p,\infty}$-norm and maximal $L^{p}$-norm, bilinear maximal $L^{p}$-norm.
\begin{proposition}\label{BLp Lp}
	Let $p \geq 2$. Suppose that $f$ has Fourier support in $B^n(0,1)$. For any $\epsilon>0$, if
	$$    \Big\|\sup_{0<t\leq r}|e^{it\Delta}f|\Big\|_{L^{p}(B_{r})} \lesssim r^{\epsilon}\|f\|_{L^{2}},  $$
	then
	$$    \|e^{it\Delta}f\|_{BL^{p,\infty}(B_{r}^{\ast})} \lesssim K^{O(1)} r^{\epsilon}\|f\|_{L^{2}} . $$
\end{proposition}

\noindent\textit{Proof.}
By the definition of $BL^{p,\infty}$-norm, we obtain
\begin{align*}
	\|e^{it\Delta}f\|_{BL^{p,\infty}(B_{r}^{\ast})}=&\Big( \sum_{B_{K}\subset B(0,r)} \max_{I_{K}\subset [0,r]} \mu_{e^{it\Delta}f}(B_{K} \times I_{K})   \Big)^{\frac{1}{p}}       \\
	\leq & \Big[ \sum_{B_{K}\subset B(0,r)} \max_{I_{K}\subset [0,r]}\Big(\sum_{\tau}\|e^{it\Delta}f_{\tau}\|^{p}_{L^{p}(B_{K}\times I_{K})}\Big)   \Big]^{\frac{1}{p}} \\
	\leq & \Big(\sum_{\tau} \sum_{B_{K}\subset B(0,r)} \max_{I_{K}\subset [0,r]}\|e^{it\Delta}f_{\tau}\|^{p}_{L^{p}(B_{K}\times I_{K})}   \Big)^{\frac{1}{p}}  \\
	\lesssim& K^{O(1)} \Big( \sum_\tau\|e^{it\Delta}f_{\tau}\|^{p}_{L^{p}_{x}L^{\infty}_{t}(B_{r}^{\ast})}\Big)^{\frac{1}{p}} \\
	\lesssim& K^{O(1)} r^{\epsilon}\|f\|_{L^{2}} .
\end{align*}%
\qed

\begin{proposition}\label{BLp BilLp}
	Let $p \geq 2$. Suppose that $f$ has Fourier support in $B^n(0,1)$.
	Let $f=\sum_{\tau} f_\tau$, where each $f_\tau$ is Fourier supported on $\tau$. For any $\epsilon>0$, if
	$$    \Big\|\sup_{0<t\leq r}\Big|e^{it\Delta}f_{\tau_{1}}e^{it\Delta}f_{\tau_{2}}\Big|^{\frac{1}{2}}\Big\|_{L^{p}(B_{r})} \lesssim r^{\epsilon}\|f\|_{L^{2}}  $$
	holds for all $\tau_1,\tau_2$ satisfying that the Fourier supports of $f_{\tau_{1}}$ and $f_{\tau_{2}}$ are separated by at least $1/K$,
	then
	$$    \|e^{it\Delta}f\|_{BL^{p,\infty}(B_{r}^{\ast})} \lesssim K^{O(1)} r^{\epsilon}\|f\|_{L^{2}} . $$
\end{proposition}

\noindent\textit{Proof.}
We recall the definition of $\mu_{e^{it\Delta}f}$:
$$ \mu_{e^{it\Delta}f}(B_{K}\times I_{K})=\min_{V_{1},...,V_{A}} \Big(  \max_{\tau \notin V_{a}} \int_{B_{K}\times I_{K}}|e^{it\Delta}f_{\tau}|^{p} dxdt  \Big).   $$
For each $\tau$, by Lemma \ref{lcp th}, there exists a positive function $\eta_K$ which is essentially supported on $B_K \times I_K$ satisfying $\|\eta_K\|_{L^1}\lesssim 1$ and
$$ |e^{it\Delta}f_{\tau}|\leq |e^{it\Delta}f_{\tau}| \ast \eta_K.  $$
Fix $B_{K} \times I_{K}$, then there exists a $\tau_{0}$ such that
$$   \max_{\tau} \int_{B_{K}\times I_{K}}|e^{it\Delta}f_{\tau}|^{p}= \int_{B_{K}\times I_{K}}|e^{it\Delta}f_{\tau_{0}}|^{p}dxdt.  $$
We choose $V'_{1},...,V'_{A}$ such that $\tau_{0} \in V'_{a}$ for some $a\in \{ 1,2,...,A \}$, then by (\ref{lcp}) and H\"older inequality,
\begin{align*}
	& \mu_{e^{it\Delta}f}(B_{K}\times I_{K})\leq  \max_{\tau \notin V'_{a}} \int_{B_{K}\times I_{K}}|e^{it\Delta}f_{\tau}|^{p}      \\
	\leq&  \max_{\tau \notin V'_{a}} \left(\int_{B_{K}\times I_{K}}|e^{it\Delta}f_{\tau}|^{p} \int_{B_{K}\times I_{K}}|e^{it\Delta}f_{\tau_0}|^{p}\right)^{\frac{1}{2}}        \\
	\lesssim& \max_{\tau \notin V'_{a}}\int_{B_{K}\times I_{K}}\left(|e^{it\Delta}f_{\tau}|\ast\eta_K\right)^{\frac{p}{2}} \left(|e^{it\Delta}f_{\tau_{0}}|\ast\eta_K\right)^{\frac{p}{2}}     \\
	\lesssim& \max_{\tau \notin V'_{a}}\int_{B_{K}\times I_{K}}\left(|e^{it\Delta}f_{\tau}|^{\frac{p}{2}} \ast\eta_K\right) \left(|e^{it\Delta}f_{\tau_{0}}|^{\frac{p}{2}} \ast\eta_K\right)    \\
	=  &  \max_{\tau \notin V'_{a}}\int \left(\int_{B_{K}\times I_{K}} \left|e^{i(t-t_1)\Delta}f_\tau(x-x_1) e^{i(t-t_2)\Delta}f_{\tau_0}(x-x_2) \right|^{\frac{p}{2}}  dxdt\right) \cdot \\
	&  \quad\quad\quad\quad\quad\quad  \eta_{K}(x_1,t_1)  \eta_{K}(x_2,t_2) dx_1dx_2dt_1dt_2  \\
	\leq&  \sum_{\substack{\tau_{1},\tau_{2} \\  d(\tau_1,\tau_2)\geq (KM)^{-1} }}
	\int \left(\int_{B_{K}\times I_{K}} \left|e^{i(t-t_1)\Delta}f_{\tau_1}(x-x_1) e^{i(t-t_2)\Delta}f_{\tau_2}(x-x_2) \right|^{\frac{p}{2}}  dxdt\right) \cdot \\
	&  \quad\quad\quad\quad\quad\quad  \eta_{K}(x_1,t_1)  \eta_{K}(x_2,t_2) dx_1dx_2dt_1dt_2  .
\end{align*}%
Define $f_{\tau_j,x_j,t_j}$ by
$$ \widehat{f_{\tau_j,x_j,t_j}} = e^{-i(t_j|\xi|^2+x_j\cdot \xi)}\widehat{f_{\tau_j}},  $$
then
$$e^{it\Delta}f_{\tau_j,x_j,t_j}=e^{i(t-t_j)\Delta}f_{\tau_j}(x-x_j).   \quad\quad j=1,~2. $$
We repeat the argument as in the proof of  Proposition \ref{BLp Lp} to obtain
\begin{align*}
	&\|e^{it\Delta}f\|_{BL^{p,\infty}(B_{r}^{\ast})}=\Big( \sum_{B_{K}\subset B(0,r)} \max_{I_{K}\subset [0,r]} \mu_{e^{it\Delta}f}(B_{K} \times I_{K})   \Big)^{\frac{1}{p}}       \\
	\leq & \Big[ \sum_{B_{K}\subset B(0,r)} \max_{I_{K}\subset [0,r]}\Big(
	\sum_{\substack{\tau_{1},\tau_{2} \\  d(\tau_1,\tau_2)\geq (KM)^{-1} }} \int\Big\|\Big|e^{it\Delta}f_{\tau_{1},x_1,t_1}e^{it\Delta}f_{\tau_{2},x_2,t_2}\Big|^{\frac{1}{2}}\Big\|^{p}_{L^{p}(B_{K}\times I_{K})}\cdot \\
	&  \quad\quad\quad\quad\quad\quad \eta_{K}(x_1,t_1)   \eta_{K}(x_2,t_2) dx_1dx_2dt_1dt_2
	\Big)   \Big]^{\frac{1}{p}}   \\
	\leq & \Big[ \sum_{\substack{\tau_{1},\tau_{2} \\  d(\tau_1,\tau_2)\geq (KM)^{-1} }}  \int \Big(\sum_{B_{K}\subset B(0,r)} \max_{I_{K}\subset [0,r]} \Big\|\Big|e^{it\Delta}f_{\tau_{1},x_1,t_1}e^{it\Delta}f_{\tau_{2},x_2,t_2}\Big|^{\frac{1}{2}}\Big\|^{p}_{L^{p}(B_{K}\times I_{K})} \Big)\cdot \\
	&  \quad\quad\quad\quad\quad\quad \eta_{K}(x_1,t_1)   \eta_{K}(x_2,t_2) dx_1dx_2dt_1dt_2
	\Big]^{\frac{1}{p}}  \\
	\lesssim& K^{O(1)} \Big( \sum_{\substack{\tau_{1},\tau_{2} \\  d(\tau_1,\tau_2)\geq (KM)^{-1} }} \int \Big\|\Big|e^{it\Delta}f_{\tau_{1},x_1,t_1}e^{it\Delta}f_{\tau_{2},x_2,t_2}\Big|^{\frac{1}{2}}\Big\|^p_{L_{x}^{p}L^{\infty}_{t}(B_{r}^{\ast})} \cdot \\
	&  \quad\quad\quad\quad\quad\quad \eta_{K}(x_1,t_1) \eta_{K}(x_2,t_2) dx_1dx_2dt_1dt_2\Big )^{\frac{1}{p}} \\
	\lesssim &K^{O(1)} r^{\epsilon}\|f\|_{L^{2}}.
\end{align*}%

\qed
\begin{remark}\label{BLp two relation}
	For $q >1/\delta$, we have
	$$ \|e^{it\Delta}f\|_{L^{p}_{x}L^{q}_{t}(B_{r}^{\ast})} \lesssim r^{\delta}\|e^{it\Delta}f\|_{L^{p}_{x}L^{\infty}_{t}(B_{r}^{\ast})}  $$
	and
	$$  \|e^{it\Delta}f\|_{BL^{p,q}(B_{r}^{\ast})} \lesssim r^{\delta}\|e^{it\Delta}f\|_{BL^{p,\infty}(B_{r}^{\ast})}. $$
	So Proposition \ref{BLp Lp} and \ref{BLp BilLp} still hold if we replace $BL^{p,\infty}$-norm with $BL^{p,q}$-norm.
\end{remark}

\noindent 4.\,\textbf{Polynomial partitioning}

Polynomial partitioning was firstly introduced by Guth \cite{G1,G2} to improve the bound of restriction conjecture in three dimension. The main idea of polynomial partitioning is using a polynomial to cut the whole space to several parts, and each part has the same contribution.
\begin{theorem}\cite[Theorem 1.4]{G1}\label{p-p orgin}
	Suppose that $F$ is a non-negative $L^1$ function on $ \mathbb{R}^{n}$. For each degree $d$, there exists a nontrivial
	polynomial $P$ of degree at most $d$ such that $\mathbb{R}^{n} \backslash Z(P)$ is a union of $ \sim_n d^{n}$ pairwise disjoint open sets $O_{i}'$(cells) and for each $i$ we have
	$$  \int_{O_{i}'} F(x)dx\sim_n d^{-n} \int F(x)dx. $$
\end{theorem}

Let $Z(P_{1},...,P_{n+1-m})$ be the set of common zeros of polynomials $P_1,\cdots, P_{m+1-n}$.
We call the variety $Z(P_{1},...,P_{n+1-m})$ to be  a transverse complete intersection if
$$  \nabla P_{1}(x) \wedge ... \wedge \nabla P_{n+1-m}(x) \neq 0,\quad \text{~~~~for ~ all~} x\in Z(P_{1},...P_{n+1-m}).   $$
If $Z=Z(P_{1},...,P_{n+1-m})$ satisfies $\deg P_i \leq d$ for each $i=1,2,...,n+1-m$, we say the degree of $Z$ is at most $d$. Du-Li in \cite{DL} gave the following polynomial partitioning theorem on the broad norm version.

\begin{theorem}\cite[Theorem 4.5]{DL}\label{p-p}
	Suppose that $f$ is a function with ${\rm supp}\widehat{f} \subset B(0,1) \subset \mathbb{R}^{n}$, $U$ is a subset of $B_{R}^{\ast}$, and $1\leq p,\lambda <\infty$. Then for each $d$, there exists a nontrivial polynomial $P$ of degree at most $O(d)$ such that $(\mathbb{R}^{n}\times \mathbb{R}) \backslash Z(P)$ is a union of $\sim d^{n+1}$ disjoint open sets $O_{i}'$ and for each $i$ we have
	$$\|e^{it\Delta}f\|^{p}_{BL^{p,\lambda}(U)} \leq C_{n} d^{n+1} \|e^{it\Delta}f\|^{p}_{BL^{p,\lambda}(U\cap O_{i}')}.  $$
	Moreover, $Z(P)$ is a finite union of transverse complete intersections.
\end{theorem}

Theorem \ref{p-p} induces one spatial decomposition $B_R^\ast= (\cup_i O_{i}') \cup Z(P)$. If we consider wave packet decomposition of $f$ on scale $R$: $f=\sum_{\theta,\nu}c_{\theta,\nu}\varphi_{\theta,\nu}$, then most of the mass of each $e^{it\Delta}\varphi_{\theta,\nu}$ is concentrated on $T_{\theta,\nu}$ in physical space. An example of $f$ says that there exists a related tube $T_{\theta,\nu}$ such that it can enter into all cells $O_i'$, which is not convenient to  carry on this divide and conquer approach. Roughly speaking, this bad situation can happen because the diameter of each cell $O_i'$ is too large. So we shrunken the scale of cells $O_i'$ to make that each tube can only enter into a small part of the shrunken cells.

Using Theorem \ref{p-p} with $\lambda=q$, we conclude that there exists a nontrivial polynomial $P$ of degree at most $O(d)$ such that $(\mathbb{R}^{n}\times \mathbb{R}) \backslash Z(P)$ is a union of $\sim d^{n+1}$ disjoint open sets $O_{i}'$ and for each $i$ we have
\begin{equation}\label{use p-p}
	\|e^{it\Delta}f\|^{p}_{BL^{p,q}(B_{R}^{\ast})} \leq C_{n} d^{n+1} \|e^{it\Delta}f\|^{p}_{BL^{p,q}(O_{i}')}.
\end{equation}

Now we define
\begin{equation}\label{signs of W and Oi}
	W:=N_{R^{\frac{1}{2}+\delta}}Z(P) \cap B_{R}^{\ast}, \quad\;~ O_{i}:=[O_{i}' \cap B_{R}^{\ast}] \backslash W,
\end{equation}
where $N_{R^{\frac{1}{2}+\delta}}Z(P)$ denotes the $R^{\frac{1}{2}+\delta}$-neighborhood of the variety $Z$. By Proposition \ref{BLp pro}, one obtains
\begin{equation}\label{Oi and W sum}
	\|e^{it\Delta}f\|^{p}_{BL^{p,q}(B_{R}^{\ast})} \leq \sum_{i} \|e^{it\Delta}f\|^{p}_{BL^{p,q}(O_{i})}+ \|e^{it\Delta}f\|^{p}_{BL^{p,q}(W)}.
\end{equation}
By the pigeonhole principle, we have at least one of the following cases holds:

\vskip0.2cm
\noindent \textbf{Cellular case}: there exist $O(d^{n+1})$ many cells $O_i$ such that for each $i$,
$$  \|e^{it\Delta}f\|^{p}_{BL^{p,q}(B_{R}^{\ast})} \lesssim d^{n+1}\|e^{it\Delta}f\|^{p}_{BL^{p,q}(O_{i})}.  $$

\noindent \textbf{Algebraic case}:
$$  \|e^{it\Delta}f\|^{p}_{BL^{p,q}(B_{R}^{\ast})} \lesssim \|e^{it\Delta}f\|^{p}_{BL^{p,q}(W)}.  $$

\noindent Therefore we can attribute the contribution on the whole region to cellular contribution and algebraic contribution. Hickman-Rogers  \cite{HR1} also introduced the polynomial partitioning theorem in lower dimensions. To build the algorithm on $BL^{p,q}$-norm in the next section, we prove relevant theorem on the version of $BL^{p,q}$-norm.

We firstly recall the definition of tangential tube introduced by Guth in \cite{G1,G2}.
\begin{definition}\label{def tang}
	Suppose that $Z=Z(P_{1},...,P_{n+1-m})$ is a transverse complete intersection in $\mathbb{R}^{n}\times \mathbb{R}$. We say that $T_{\theta,\nu}$ is concentrated on scale $r$ wave packets which are $r^{-\frac{1}{2} + \delta_{m}}$-tangent to $Z$ in $B_{r}$ if the following two conditions hold:
	
	\begin{itemize}
		\item[{\rm(i)}] {\bf Distance condition}:
		$$ T_{\theta,\nu} \subset N_{r^{\frac{1}{2}+ \delta_{m}}}Z\cap B_{r}.$$
		\item[{\rm(ii)}] {\bf Angle condition}: for any $x\in T_{\theta,\nu}$ and $z \in Z\cap B_{r} $ with $|z-x| \lesssim r^{1/2+\delta_{m}}$, one has
		$$ \angle(G(\theta),T_{z}Z) \lesssim r^{-\frac{1}{2} + \delta_{m}}.$$
	\end{itemize}
	Let
	\begin{align*} \mathbb{T}_Z[r]:=\Big\{  (\theta,\nu)~\big|~&T_{\theta,\nu} \text{~is concentrated on scale $r$ wave packets}\\
		&\text{ which are  $r^{-\frac{1}{2} + \delta_{m}}$ tangent to}~ Z \Big\}.\end{align*}
	We say that $f$ is concentrated in wave packets from $ \mathbb{T}_Z[r]$ if
	$$ \sum_{(\theta,\nu) \notin\mathbb{T}_Z[r] } \|f_{\theta,\nu}\|_{L^2} \leq O(r^{-N})\|f\|_{L^2},  \;\;\; \forall~ N>0.  $$
\end{definition}

\begin{theorem}\label{p-p low}
	Fix $r \gg 1$, $d\in \mathbb{N}$ and $1\leq p,\lambda <\infty$. Suppose that $Z$ is an $m$-dimensional transverse complete intersection of degree at most $d$. Suppose that $f$ satisfies  ${\rm supp}\widehat{f} \subset B(0,1) \subset \mathbb{R}^{n}$, and $f$ is concentrated in wave packets from $ \mathbb{T}_Z[r]$. Then there exists a constant $D=D_Z=D(\epsilon,d)$ such that at least one of the following cases holds:
	
	\noindent \textbf{Cellular case}: there exist $O(D^m)$ shrunken cells $O_i$ in $\mathbb{R}^{n} \times \mathbb{R}$ such that for each $i$,
	$$  \|e^{it\Delta}f\|^{p}_{BL^{p,\lambda}(B_{r})} \lesssim  D^{m}\|e^{it\Delta}f\|^{p}_{BL^{p,\lambda}(O_{i})}.  $$

	\noindent \textbf{Algebraic case}: there exists an $(m-1)$-dimensional transverse complete intersection $Y$ of degree at most $O(D)$ such that
	$$  \|e^{it\Delta}f\|^{p}_{BL^{p,\lambda}(B_{r} \cap N_{r^{1/2+\delta_{m}}}Z)} \lesssim \|e^{it\Delta}f\|^{p}_{BL^{p,\lambda}(B_{r} \cap N_{r^{1/2+\delta_{m}}}Y)}.  $$
\end{theorem}

\noindent\textit{Proof.}
Suppose that the algebraic case does not occur. We say that a ball $B=B(x_0,r^{1/2+\delta_m})\subset N_{r^{1/2+\delta_{m}}} Z\cap B_r $ is regular if, on each connected component of $Z\cap B(x_0,r^{1/2+\delta_m})$, the tangent space $TZ$ is constant up to angle $1/100$. By the hypothesis and Guth's argument such as in \cite[Section 8]{G2}, we have the regular balls contain most of the mass of $e^{it\Delta}f$, i.e.
$$   \|e^{it\Delta}f\|_{BL^{p,\lambda}(B_{r})} \lesssim \|e^{it\Delta}f\|_{BL^{p,\lambda}(\cup_{\text{regular}~B}B)}.  $$
For each regular ball $B$, we pick a point $z\in Z\cap B_{r^{1/2+\delta_m}}$ and define $V_B$ to be the $m$-dimensional tangent plane $T_z Z$. For each $m$-dimensional tangent plane $V$, we use $\mathfrak{B}_V$ to denote the set of regular balls so that $\angle(V_B,V)\leq 1/100$. By the pigeonhole principle, there exists a plane $V$ such that
$$   \|e^{it\Delta}f\|_{BL^{p,\lambda}(B_{r})} \lesssim \|e^{it\Delta}f\|_{BL^{p,\lambda}(\cup_{B\in \mathfrak{B}_V}B)}.  $$
Set $N_1=\cup_{B\in \mathfrak{B}_V}B$. Therefore, it suffices to consider the contribution on $N_1$.

Since all wave packets of $e^{it\Delta}f$ on scale $r$ are contained in the $\thickapprox r^{1/2}$-neighborhood of $V$, $|e^{it\Delta}f(x)|$ is essentially constant along a certain direction which is roughly normal to $V$. Let $P_V: V\rightarrow \mathbb{R}$ denote a polynomial defined on $V$, and $\pi: \mathbb{R}^{n+1}\rightarrow V$ be the orthogonal projection. Then we can extend $V$ to a polynomial $P$ on $\mathbb{R}^{n+1}$ by setting $P(x):=P_V(\pi(x))$. We denote the collection of such polynomials by $\mathbb{P}_V$. Repeating the proof of Theorem 4.5 in \cite{DL}, we can find a polynomial $P \in \mathbb{P}_V$ of degree at most $D=D(\epsilon,d)$ such that $N_1=(\cup_i O_i') \cup Z(P)$ satisfying $\#O'_i \sim D^m$, and for each $i$,
\begin{equation}\label{lllll}
	\|e^{it\Delta}f\|^{p}_{BL^{p,\lambda}(B_{r})} \sim \|e^{it\Delta}f\|^{p}_{BL^{p,\lambda}(N_1)} \lesssim  D^{m}\|e^{it\Delta}f\|^{p}_{BL^{p,\lambda}(O'_{i})}.
\end{equation}
Define
$$     W:=N_{r^{\frac{1}{2}+\delta_m}}Z(P) \cap B_{r}, \quad\;~ O_{i}:=[O_{i}' \cap B_{r}] \backslash W.   $$
By Proposition {\ref{BLp pro}}, we obtain
$$  \|e^{it\Delta}f\|^{p}_{BL^{p,\lambda}(B_{r})} \leq \sum_{i}\|e^{it\Delta}f\|^{p}_{BL^{p,\lambda}(O_{i})} +\|e^{it\Delta}f\|^{p}_{BL^{p,\lambda}(W)}. $$
Since the algebraic case does not occur, the contribution from $W$ is negligible. By the pigeonhole principle and \eqref{lllll}, there exist $O(D^m)$ cells $O_i$ such that for each $i$,
$$  \|e^{it\Delta}f\|^{p}_{BL^{p,\lambda}(B_{r})} \lesssim  D^{m}\|e^{it\Delta}f\|^{p}_{BL^{p,\lambda}(O_{i})}.  $$

\qed

\noindent 5.\,\textbf{The vanishing property of the broad norm}

In the final part of this section, we prove one important property of the broad norm named the  vanishing property. For the $BL_{k,A}^p$-norm, Guth \cite{G2} has proved this property. Now we show this property also holds for general $BL_{k,A}^{p,q}$-norm.

\begin{theorem}\label{vp}
Let $r \gg 1$, $\delta \ll\epsilon$ and $1 \leq m <k \leq n$, and let $Z$ be an $m$-dimensional transverse complete intersection. Suppose that $f$ is concentrated on scale $r$ wave packets which are $r^{-\frac{1}{2}+\delta}$-tangent to $Z$ on $B_r$, then
$$    \|e^{it\Delta}f\|_{BL_{k,A}^{p,q}(B_r)} =O(r^{-N})\|f\|_{L^2}   $$
for any $N>0$.
\end{theorem}

\noindent \textit{Proof.} We consider wave packet decomposition of $f$ on scale $r$: $f=\sum_{\theta,\nu}c_{\theta,\nu}\varphi_{\theta,\nu}$.  Since  $f$ is concentrated in wave packets from $ \mathbb{T}_Z[r]$, it suffices to consider $(\theta,\nu) \in  \mathbb{T}_Z[r]$. For each $(\theta,\nu) \in  \mathbb{T}_Z[r]$,  we can choose $z_{0} \in Z\cap B_r \cap N_{O(r^{1/2+ \delta})}T_{\theta,\nu}$ such that
$$ \angle (G(\theta),T_z Z) \lesssim r^{-\frac{1}{2} + \delta}.$$
This fact implies for any $\tau \supset \theta$,
$$ \angle (G(\tau),T_z Z) \leq K^{-1}.$$
Since $T_z Z$ is $m$-dimensional, note $m <k$, we have $T_z Z$ is also a ($k-1$)-dimensional subspace. So such $\tau$'s do not contribute to $\mu_{e^{it\Delta}f}(B_{K}\times I_{K})$, and then
$$    \|e^{it\Delta}f\|_{BL_{k,A}^{p,q}(B_r)} =O(r^{-N})\|f\|_{L^2}   $$
for any $N>0$.

\qed

\section{Two algorithms}\label{section4}

From this section, we start to prove Theorem \ref{th4}. We will adopt two algorithms introduced by Hickman-Rogers in \cite{HR1}.
The method developed by Guth \cite{G2} in  the study of  restriction conjecture relies on polynomial partitioning and scale induction.
Hickman-Rogers used the algorithm instead of scale induction to give a new proof of Guth's result.
The main advantage of the algorithm is to lead to a more detailed geometry analysis on different scales.

Now we introduce the first algorithm. It is a dimensional reduction, essentially passing from an $m$-dimensional to an $(m-1)$-dimensional situation. More precisely, we start with a function $f$ which is concentrated in wave packets from $ \mathbb{T}_Z[r]$,  and $Z$ is an $m$-dimensional transverse complete intersection. Using Theorem \ref{p-p low}, we can divide  $\|e^{it\Delta}f\|_{BL^{p,q}(B_{r})}$ into the cellular and algebraic case. The algebraic case can further separated into the transverse and tangential case according to the angle condition. If it's the tangential case, we already pass from an $m$-dimensional to an $(m-1)$-dimensional situation, then the algorithm stops. Otherwise, we repeat the above steps to each cells produced by the cellular and transverse cases.
\vskip0.5cm
\noindent \textbf{[Alg-1]}-{The first algorithm}:\;  Let $p \geq 2$,\ $0 < \epsilon  \ll 1$, and $\delta,\delta_{n},...,\delta_1,\delta_{0}$ satisfy
$$  \epsilon^{C} \leq \delta \ll \delta_{n} \ll \delta_{n-1} \ll ...\ll \delta_{1} \ll \delta_{0}   \ll \epsilon     $$
for some constant $C$.
\vskip0.5cm

\noindent \underline{\bf Input}.
\begin{itemize}
	\item  A ball $B_{r}\subset \mathbb{R}^{n+1}$ with radius $r \gg 1.$
	\item  A transverse complete intersection\ $Z$ of dimension $m \geq 2.$
	\item  A function $f$ satisfying supp$\widehat{f} \subset B^n(0,1)$, and $f$ is concentrated on scale $r$ wave packets which are\ $r^{-1/2+\delta_{m}}$-tangent to\ $Z$ in $B_{r}$.
\end{itemize}
\vskip0.5cm

\noindent \underline{\bf Output}. The $j$-th step of recursion will produce:

\begin{itemize}
	\item A word $   \mathcal{E}_j $ of length $j$ in the alphabet \{{\bf a,c}\}. Here {\bf a} is an abbreviation of ``algebraic" and {\bf c}  ``cellular". The algorithm is realized by repeated application of Theorem \ref{p-p low}. More precisely, suppose that the algorithm has ran through $(j-1)$ steps, and produced a word $\mathcal{E}_{j-1}$ of length $j-1$ and a family of subsets $\mathcal{O}_{j-1}$ (see its definition below), the we again use Theorem \ref{p-p low} on each $O_{j-1} \in \mathcal{O}_{j-1}$. If more than 1/2 fraction of the $\mathcal{O}_{j-1}$ belong to the cellular case, then we define $\mathcal{E}_j$ as a word of length $j$, where the first $(j-1)$ letters are the same as $\mathcal{E}_{j-1}$ and the $j$-th letter is {\bf c}. Otherwise, we defined the $j$-th letter of $\mathcal{E}_j$ to be {\bf a}.
	
	\item A spatial scale $\rho_{j} \geq 1$. Let $\tilde{\delta}_{m-1}$ satisfy
	$$   (1-\tilde{\delta}_{m-1})\Big(\frac{1}{2}+\delta_{m-1}\Big)=\frac{1}{2}+\delta_{m} ,$$
	and $\sigma_{k}: [0,\infty) \rightarrow [0,\infty)$ satisfy
	\begin{align*}
		\sigma_{k}(\rho):=
		\begin{cases}
			\frac{\rho}{2},~~~~~\ \ \ \ \ \ \ \ \ \ \ \ \ \  \text{if the } k \text{-th letter of~} \mathcal{E}_{j} \text{~is~} {\bf c},  \\
			\rho^{1-\tilde{\delta}_{m-1}},~~~~~ \ \  \ \ \ \  \text{if the } k \text{-th letter of~} \mathcal{E}_{j} \text{~is~} {\bf a},
		\end{cases}
	\end{align*}
	for each $1\leq k\leq j$. With these definitions, take
	$$    \rho_{j}:=\sigma_{j}\circ...\circ\sigma_{1}(r).   $$
	Then $\rho_{j}$ satisfies:
	$$   \rho_{j} \leq r^{(1-\tilde{\delta}_{m-1})^{\#{\bf a}(j)}} \;\;\text{and}\;\;\rho_{j}\leq \frac{r}{2^{\#{\bf c}(j)}},   $$
	where\ $\#{\bf a}(j)$ and\ $\#{\bf c}(j)$ denote the number of occurrences of {\bf a} and {\bf c} in the $\mathcal{E}_{j}$.
	\vskip0.2cm
	
	\item A family of subsets $\mathcal{O}_{j}$ of $\mathbb{R}^{n+1}$ which will be referred to as cells. Each cell $O_{j} \in \mathcal{O}_{j}$ has diameter at most\ $\rho_{j}$.
	
	\item A collection of functions $(f_{O_{j}})_{O_{j} \in \mathcal{O}_{j}}$. Each $f_{O_{j}}$ is concentrated on scale $\rho_{j}$ wave packets which are\ $\rho_{j}^{-1/2+\delta_{m}}$-tangent to some translation of $Z$ on (a ball of radius $\rho_{j}$ containing) $O_{j}$.
	
	\item A large integer $d\in \mathbb{N}$ which depends only on the admissible parameters and $\deg Z$.
\end{itemize}
\vskip0.15cm

\noindent Then the following properties hold:
\vskip0.15cm

\noindent \textbf{Property I}. Most of the mass of $\|e^{it\Delta}f\|^{p}_{BL^{p,q}(B_{r})}$ is concentrated on the $O_{j} \in \mathcal{O}_{j}$:
\begin{equation}
	\|e^{it\Delta}f\|^{p}_{BL^{p,q}(B_{r})} \leq  C^{\text{I}}_{j,\delta}(d,r) \sum_{O_{j} \in \mathcal{O}_{j}}    \|e^{it\Delta}f_{O_{j}}\|^{p}_{BL^{p,q}(O_{j})}+ O(r^{-N})\|f\|_{L^2}.
\end{equation}

\noindent \textbf{Property II}. The functions $f_{O_{j}}$ satisfy
\begin{equation}
	\sum_{O_{j} \in \mathcal{O}_{j}} \|f_{O_{j}}\|^{2}_{L^{2}} \leq C^{\text{II}}_{j,\delta}(d,r)  d^{\#{\bf c}(j)}\|f\|^{2}_{L^{2}}.
\end{equation}

\noindent \textbf{Property III}. Each $f_{O_{j}}$ satisfies
\begin{equation}
	\|f_{O_{j}}\|^{2}_{L^{2}} \leq C^{\text{III}}_{j,\delta}(d,r)\Big(\frac{r}{\rho_{j}}\Big)^{-\frac{n+1-m}{2}}d^{-\#{\bf c}(j)(m-1)}\|f\|^{2}_{L^{2}}.
\end{equation}
Here
\begin{equation*}\left\{\begin{aligned}
		&   C^{\text{I}}_{j,\delta}(d,r):=d^{\#{\bf c}(j) \delta}(\log r)^{2\#{\bf a}(j)(1+\delta)},\\
		& C^{\text{II}}_{j,\delta}(d,r):=d^{\#{\bf c}(j) \delta}\text{Poly}(d)^{\#{\bf a}(j)(1+\delta)}, \\
		& C^{\text{III}}_{j,\delta}(d,r):=d^{j \delta}r^{\overline{C}\#{\bf a}(j)\delta_{m}}, \end{aligned}\right.\end{equation*}
where $\text{Poly}(d):=d^{n+1}$ and $\overline{C}$ is some suitably chosen large constant. One easily verifies that
\begin{equation}\label{sign 1,2,3}
	C^{\text{I}}_{j,\delta}(d,r),\; \;C^{\text{II}}_{j,\delta}(d,r),\; \;C^{\text{III}}_{j,\delta}(d,r) \lesssim_{d,\delta}r^{\delta_{0}}d^{\#{\bf c}(j)\delta}.
\end{equation}
\vskip0.25cm

\noindent \underline{\bf The first step}. The algorithm shall start by taking:
\begin{itemize}
	\item $\mathcal{E}_0:=\varnothing$ to be the empty word;
	\item $\rho_0:=r$;
	\item $\mathcal{O}_{0}:=\{ O_0\}$, where $O_0:=N_{r^{1/2+\delta_m}}Z \cap B_r$;
	\item $f_{O_0}:=f$.
	\item A large integer $d\in \mathbb{N}$ to be determined later.
\end{itemize}
Under this setting, one easily verifies the validity of {\bf Property I, II}, and {\bf III}.
In fact, {\bf Property I} holds due to the hypothesis of $f$. {\bf Property II} and {\bf III}  hold trivially.
\vskip0.25cm

\noindent \underline{\bf The $(j+1)$-th step}. Let $j \geq 1$. We assume that the recursive algorithm has ran through $j$ steps. Since each function $f_{O_j}$ is concentrated on scale $\rho_j$ wave packets which are\ $\rho_j^{-1/2+\delta_{m}}$-tangent to some translation of $Z$ in $B_{\rho_j}$, we again apply Theorem \ref{p-p low} to bound $\|e^{it\Delta}f_{O_{j}}\|^{p}_{BL^{p,q}(O_{j})}$ for each $O_{j} \in \mathcal{O}_{j}$. One of two things can happen: either \textbf{[Alg-1]} terminates if the diameter of each cell is very small or it terminates if most of the mass of $\|e^{it\Delta}f_{O_{j}}\|^{p}_{BL^{p,q}(O_{j})}$ is concentrated on some $(m-1)$-dimensional transverse complete intersections. We give the exact stopping conditions as follows.
\vskip0.25cm

\noindent  \underline{\bf Stopping conditions}. The algorithm will be stopped if one of both conditions occurs:
\begin{itemize}
	\item[\textbf{[tiny]}]\;  The algorithm terminates if $\rho_{j} \leq r^{\tilde{\delta}_{m-1}}$ .
	
	\item[ \textbf{[tang]}]\;  The algorithm terminates if
	\begin{equation}\label{tang 1 condition}
		\sum_{O_{j} \in \mathcal{O}_{j}} \|e^{it\Delta}f_{O_{j}}\|^{p}_{BL^{p,q}(O_{j})} \leq C_{\mathrm{tang}} \sum_{S \in \mathcal{S}} \|e^{it\Delta}f_{S}\|^{p}_{BL^{p,q}(B_{\tilde{\rho}})}
	\end{equation}
	and
	\begin{equation}
		\sum_{S \in \mathcal{S}} \|f_{S}\|^{2}_{L^{2}} \leq C_{\mathrm{tang}} r^{(n+1)\tilde{\delta}_{m-1}} \sum_{O_{j} \in \mathcal{O}_{j}}\|f_{O_{j}}\|^{2}_{L^{2}}
	\end{equation}
	hold for $ \tilde{\rho}:=\rho_{j}^{1-\tilde{\delta}_{m-1}}$ and for some constant $C_{\mathrm{tang}}$. Here $\mathcal{S}$ is a collection of $(m-1)-$dimensional transverse complete intersections in $\mathbb{R}^{n+1}$ all of degree\ $O(d)$, and $f_{S}$ is a function to each $S \in \mathcal{S}$ which is $\tilde{\rho}^{-1/2+\delta_{m-1}}$-tangent to\ $S$ on $B_{\tilde{\rho}}$.
\end{itemize}

\begin{remark}
	\noindent{\rm (a)}\; We can use Theorem \ref{p-p low} and scale induction to reformulate this algorithm. For details one can refer to Hickman-Rogers
	\cite[Section 9]{HR1}.
	
	\noindent{\rm (b)} There exists a additional local estimate in {\bf Property III} in {\rm [\bf Alg-1]} of Hickman-Rogers. The main effect of such estimate is to establish  $L^\infty$-estimate from $L^2$-estimate at the smallest possible cost due to the hypothesis of restriction conjecture. In this paper, we always assure $f\in L^2$, so the local estimate is not essential in our argument though it also holds in this case.
\end{remark}
Next we introduce the second algorithm. It consists of repeated application of the first algorithm to reduce to an appropriate dimensional case.
\vskip0.25cm

\noindent \textbf{[Alg-2]}-The second algorithm:\; Let $\{p_{\ell}\}_{\ell=2}^{n+1}$ denote Lebesgue exponents such that
$$    p_{2} \geq p_{3} \geq ...\geq p_{n+1}:=p\geq 2.  $$
Let \ $0 \leq \alpha_{\ell},\beta_{\ell} \leq 1$ be the numbers defined in terms of $p_\ell$ by
$$      \frac{1}{p_{\ell}}:=\frac{1-\alpha_{\ell-1}}{2} +\frac{\alpha_{\ell-1}}{p_{\ell-1}},\quad~ \beta_{\ell}:= \prod_{i=\ell}^{n}\alpha_{i}~\;\;~\text{for} \;\;
3\leq \ell \leq n, $$
and $\alpha_{n+1}:=\beta_{n+1}:=1$.
\vskip0.25cm

\noindent \underline{\bf Input}.  Fix\ $R \gg 1$. Let\ $f$ satisfy supp$\widehat{f} \subset B^n(0,1)$ and the non-degeneracy hypothesis
\begin{equation}\label{ndh}
\|e^{it\Delta}f\|_{BL^{p,q}(B^{\ast}_{R})} \geq CR^{\epsilon}\|f\|_{L^{2}}.
\end{equation}
\vskip0.25cm

\noindent \underline{\bf Output}. The\ $(n+2-\ell)$-th step of the recursion will produce:

\begin{itemize}
	\item Scales\ $\vec{r}_{\ell}=(r_{n+1},...,r_{\ell})$, satisfying\ $R = r_{n+1} > r_{n} > ...> r_{\ell}.$ Large and non-admissible parameters $\vec{D}_{\ell}=(D_{n+1},...,D_{\ell}).$ Here $D_i$, $\ell \leq i \leq n+1$, are defined by the upper bound of (\ref{sign 1,2,3}).
	\item A family transverse complete intersections $\mathcal{S}_{\ell}$. Each $S_{i} \in \mathcal{S}_{i}$, $\ell \leq i\leq n+1$, satisfies $\deg S_{i}=i$.
	\item  A function $f_{S_{\ell} }$ associated with each $S_{\ell} \in \mathcal{S}_{\ell}$ is concentrated on scale\ $r_{\ell}$ wave packets which are $r_{\ell}^{-1/2+\delta_{\ell}}$-tangent to  $S_{\ell}$ in\ $B_{r_{\ell}}$.
\end{itemize}
\vskip0.15cm

Then the following properties hold:
\vskip0.20cm
\noindent \textbf{Property 1}.
\begin{equation}
	\|e^{it\Delta}f\|_{BL^{p,q}(B^{\ast}_{R})} \lessapprox M(\vec{r}_{\ell},\vec{D}_{\ell}) \|f\|^{1-\beta_{\ell}}_{L^{2}}\Big( \sum_{S_{\ell} \in \mathcal{S}_{\ell}}  \|e^{it\Delta}f_{S_{\ell} }\|^{p_{\ell}}_{BL^{p_{\ell},q}(B_{r_{\ell}})}              \Big)^{\frac{\beta_{\ell}}{p_{\ell}}}  ,
\end{equation}
where
$$  M(\vec{r}_{\ell},\vec{D}_{\ell}):= \Big(\prod_{i=l}^{n} D_{i}\Big)^{(n+1-\ell)\delta}\Big( \prod_{i=l}^{n}r_{i}^{\frac{i}{2(i+1)}(\beta_{i+1}-\beta_{i})}  D_{i}^{\frac{1}{2}(\beta_{i+1}-\beta_{\ell})}\Big).    $$

\noindent \textbf{Property 2}. For\ $\ell \leq n$,
\begin{equation}
	\sum_{S_{\ell} \in \mathcal{S}_{\ell}} \|f_{S_{\ell}}\|^{2}_{L^{2}} \lessapprox  D^{1+\delta}_{\ell} \sum_{S_{\ell+1} \in \mathcal{S}_{\ell+1}}  \|f_{S_{\ell+1}}\|^{2}_{L^{2}}.
\end{equation}

\noindent \textbf{Property 3}. For\ $\ell \leq n$,
\begin{equation}
	\max_{S_{\ell} \in \mathcal{S}_{\ell}  } \|f_{S_{\ell}}\|^{2}_{L^{2}} \lessapprox \Big( \frac{r_{\ell+1}}{r_{\ell}}\Big)^{-\frac{n-\ell}{2}}D^{-\ell+\delta}_{\ell}\max_{S_{\ell+1} \in \mathcal{S}_{\ell+1}  }  \|f_{S_{\ell+1}}\|^{2}_{L^{2}}.
\end{equation}
\vskip0.5cm

\noindent \underline{\bf The first step}. The algorithm shall start by taking:
\begin{itemize}
	\item $r_{n+1}:=R$ and $D_{n+1}:=1$;
	\item $\mathcal{S}_{n+1}:=\{ S_{n+1}\}$, where $S_{n+1}:= \mathbb{R}^{n+1}$;
	\item $\mathcal{O}_{0}:=\{ O_0\}$, where $O_0:=N_{r^{1/2+\delta_m}}Z \cap B_r$;
	\item $f_{S_{n+1}}:=f$ is a function satisfying the non-degeneracy hypothesis, and $f$ is concentrated on scale $R$ wave packets which are $R^{-1/2+\delta_{n+1}}$-tangent to the $(n+1)$-dimensional variety $\mathbb{R}^{n+1}$ in $B^\ast_R$.
\end{itemize}

With these definitions, {\bf Property 1, 2} and {\bf 3} hold trivially.
\vskip0.5cm
\noindent \underline{\bf The $(n+3-\ell)$-th step}. Let $\ell \geq 1$. We assume that the algorithm has ran through $(n+2-\ell)$ steps. Since  each function $f_{S_\ell}$ is concentrated on scale $r_{\ell}$ wave packets which are\ $r_\ell^{-1/2+\delta_{\ell}}$-tangent to $S_\ell$ in $B_{r_\ell}$, we repeat to apply \textbf{[Alg-1]} to bound $\|e^{it\Delta}f_{S_{\ell}}\|^{p_\ell}_{BL^{p_\ell,q}(B_{r_\ell})}$ for each $S_{\ell} \in \mathcal{S}_{\ell}$. One of two things can happen: either \textbf{[Alg-2]} terminates if most of $S_{\ell} \in \mathcal{S}_{\ell}$ ended with the stopping condition \textbf{[tiny]} after using \textbf{[Alg-1]} or it terminates due to the stopping condition \textbf{[tang]}. This recursive process terminates if the contributions from terms of the former type dominate. We give the exact stopping conditions as follows.

\vskip0.5cm

\noindent  \underline{\bf Stopping conditions}. The algorithm will be stopped if the following condition occurs:
\begin{itemize}
	\item[\textbf{[tiny-dom]}]\; Suppose that the inequality
	\begin{equation}
		\sum_{S_{\ell} \in \mathcal{S}_{\ell}}  \|e^{it\Delta}f_{S_{\ell} }\|^{p_{\ell}}_{BL^{p_{\ell},q}(B_{r_{\ell}})}  \leq 2\sum_{S_{\ell} \in \mathcal{S}_{\ell,\mathrm{tiny}}}  \|e^{it\Delta}f_{S_{\ell} }\|^{p_{\ell}}_{BL^{p_{\ell},q}(B_{r_{\ell}})}
	\end{equation}
	holds, where the right-hand summation is restricted to those $S_{\ell} \in \mathcal{S}_{\ell}$ for which \textbf{[Alg-1]} terminates owing to the stopping condition \textbf{[tiny]}. Then \textbf{[Alg-2]} terminates.
\end{itemize}

{\bf Property 2} and {\bf 3} of \textbf{[Alg-2]} can be obtained by repeated applications of \textbf{Property $\text{II}$} and \textbf{$\text{III}$} in \textbf{[Alg-1]}.
We give the proof of {\bf Property 1}, which relies on the following result. This result can be proved by the same argument as in \cite{DGLZ}, so here we omit the proof.

\begin{theorem}\label{Du Zhang result}
	Let $n\geq 2$, $\delta \ll \epsilon$, and $k $ be a dimension in the range $2 \leq k \leq n$. Suppose that $Z=Z(P_1,...P_{n+1-k})$ is a transverse complete intersection where $\deg (P_i) \leq D_Z=R^{\delta_{\deg} }$. Here $\delta_{\deg} \ll \delta$ is a small parameter. Suppose that $f$ has Fourier support in $B^n(0,1)$, and $f$ is concentrated on scale $R$ wave packets which are $R^{-\frac{1}{2}+\delta}$-tangent to $Z$ on $B_R$. Then we have
	\begin{equation}
		\Big\|\sup_{0<t\leq R} |e^{it\Delta}f|  \Big\|_{L^{2}(B^n(0,R))} \lessapprox R^{\frac{k}{2(k+1)}}\|f\|_{L^{2}}
	\end{equation}
	for all $R \geq 1$.
\end{theorem}

\noindent \textit{The proof of Property 1.}
If the stopping condition \textbf{[tiny-dom]} not happened, then
\begin{equation}
	\sum_{S_{\ell} \in \mathcal{S}_{\ell}}  \|e^{it\Delta}f_{S_{\ell} }\|^{p_{\ell}}_{BL^{p_{\ell},q}(B_{r_{\ell}})}  \leq 2\sum_{S_{\ell} \in \mathcal{S}_{\ell,\mathrm{tang}}}  \|e^{it\Delta}f_{S_{\ell} }\|^{p_{\ell}}_{BL^{p_{\ell},q}(B_{r_{\ell}})}
\end{equation}
holds, where the right-hand summation is restricted to those $S_{\ell} \in \mathcal{S}_{\ell}$ for which \textbf{[Alg-1]} terminates owing to the stopping condition \textbf{[tang]}.
For each $S_{\ell} \in \mathcal{S}_{\ell,\mathrm{tang}}$, let $\mathcal{S}_{\ell-1}[S_{\ell}]$ denote a collection of transverse complete intersections in $\mathbb{R}^{n+1}$ of dimension $\ell-1$ which is produced by the stopping condition \textbf{[tang]} in \textbf{[Alg-1]} on $\|e^{it\Delta}f_{S_{\ell} }\|^{p_{\ell}}_{BL^{p_{\ell},q}(B_{r_{\ell}})}$, then by {\bf Proposition I} in \textbf{[Alg-1]},
$$ \|e^{it\Delta}f_{S_{\ell} }\|^{p_{\ell}}_{BL^{p_{\ell},q}(B_{r_{\ell}})} \lessapprox D_{\ell-1}^{\delta} \sum_{S_{\ell-1} \in \mathcal{S}_{\ell-1}[S_{\ell}]}  \|e^{it\Delta}f_{S_{\ell-1} }\|^{p_{\ell}}_{BL^{p_{\ell},q}(B_{r_{\ell-1}})}.   $$
Define
$$  \mathcal{S}_{\ell-1}:=\{ S_{\ell-1} : S_{\ell} \in \mathcal{S}_{\ell,\mathrm{tang}} \text{~and~} S_{\ell-1} \in \mathcal{S}_{\ell-1}[S_{\ell}]\}. $$
By induction on $\ell$, we have
$$  \|e^{it\Delta}f\|_{BL^{p,q}(B_{R}^{\ast})} \lessapprox D_{\ell-1}^{\delta}M(\vec{r}_{\ell},\vec{D}_{\ell})\|f\|_{L^{2}}^{1-\beta_{\ell}} \Big( \sum_{S_{\ell-1} \in \mathcal{S}_{\ell-1}}  \|e^{it\Delta}f_{S_{\ell-1} }\|^{p_{\ell}}_{BL^{p_{\ell},q}(B_{r_{\ell-1}})}  \Big)^{\frac{\beta_{\ell}}{p_{\ell}}}.  $$
On the other hand, by Proposition \ref{BLp pro} and H\"older's inequality, we have
\begin{align*}
	&\Big( \sum_{S_{\ell-1} \in \mathcal{S}_{\ell-1}}  \|e^{it\Delta}f_{S_{\ell-1} }\|^{p_{\ell}}_{BL^{p_{\ell},q}(B_{r_{\ell-1}})}  \Big)^{\frac{1}{p_{\ell}}} =
	\Big\| \|e^{it\Delta}f_{S_{\ell-1} }\|_{BL^{p_{\ell},q}(B_{r_{\ell-1}})}  \Big\|_{l^{p_{\ell}}(\mathcal{S}_{\ell-1})}    \\
	\leq&   \Big\| \|e^{it\Delta}f_{S_{\ell-1} }\|^{1-\alpha_{\ell-1}}_{BL^{2,q}(B_{r_{\ell-1}})}  \|e^{it\Delta}f_{S_{\ell-1} }\|^{\alpha_{\ell-1}}_{BL^{p_{\ell-1},q}(B_{r_{\ell-1}})}   \Big\|_{l^{p_{\ell}}(\mathcal{S}_{\ell-1})}     \\
	\leq&    \Big\| \|e^{it\Delta}f_{S_{\ell-1} }\|_{BL^{2,q}(B_{r_{\ell-1}})}  \Big\|^{1-\alpha_{\ell-1}}_{l^{2}(\mathcal{S}_{\ell-1})}    \Big\| \|e^{it\Delta}f_{S_{\ell-1} }\|_{BL^{p_{\ell-1},q}(B_{r_{\ell-1}})}   \Big\|^{\alpha_{\ell-1}}_{l^{p_{\ell-1}}(\mathcal{S}_{\ell-1})} .
\end{align*}%
From Theorem \ref{Du Zhang result}  it follows
$$ \|e^{it\Delta}f_{S_{\ell-1} }\|_{L_{x}^{2}L^{\infty}_{t}(B_{r_{\ell-1}})} \lessapprox r_{\ell-1}^{\frac{\ell-1}{2\ell}}\|f_{S_{\ell-1}}\|_{L^{2}},    $$
and so we have  by Proposition \ref{BLp Lp} and Remark \ref{BLp two relation}
$$ \|e^{it\Delta}f_{S_{\ell-1} }\|_{BL^{2,q}(B_{r_{\ell-1}})} \lessapprox  \|e^{it\Delta}f_{S_{\ell-1} }\|_{BL^{2,\infty}(B_{r_{\ell-1}})}  \lessapprox r_{\ell-1}^{\frac{\ell-1}{2\ell}}\|f_{S_{\ell-1}}\|_{L^{2}}. $$
This inequality implies
$$    \Big\| \|e^{it\Delta}f_{S_{\ell-1} }\|_{BL^{2,q}(B_{r_{\ell-1}})}  \Big\|_{l^{2}(\mathcal{S}_{\ell-1})}  \lessapprox  r_{\ell-1}^{\frac{\ell-1}{2\ell}} \Big(  \prod_{i=l-1}^{n} D_{i}^{1+\delta} \Big)^{\frac{1}{2}}\|f\|_{L^{2}}.  $$
Combining all estimates, one concludes
$$  \|e^{it\Delta}f\|_{BL^{p,q}(B_{R}^{\ast})} \lessapprox M(\vec{r}_{\ell-1},\vec{D}_{\ell-1})\|f\|_{L^{2}}^{1-\beta_{\ell-1}} \Big( \sum_{S_{\ell-1} \in \mathcal{S}_{\ell-1}}  \|e^{it\Delta}f_{S_{\ell-1} }\|^{p_{\ell-1}}_{BL^{p_{\ell-1},q}(B_{r_{\ell-1}})}  \Big)^{\frac{\beta_{\ell-1}}{p_{\ell-1}}}.  $$

\qed

\section{The final stage}\label{section5}

In this section, we apply the algorithms in the previous section to prove Theorem \ref{th4}. We assume that \textbf{[Alg-2]} terminates the step $m$, then $m\geq 2$. In fact, if the recursive process terminates the step $m=1$, then each function $f_{S_1}$ is concentrated on wave packets which are tangent to some transverse complete intersection of dimension 1. By Theorem \ref{vp}, we have
$$    \|e^{it\Delta}f_{S_1}\|_{BL^{p,q}(B_{r_1})} =O(r_1^{-N})\|f_{S_1}\|_{L^2}   $$
for any $N>0$. Then it easily follows that
$$      \|e^{it\Delta}f\|_{BL^{p,q}(B_{R})} =O(R^{-N})\|f\|_{L^2}, $$
which contradicts the non-degeneracy hypothesis (\ref{ndh}).

Recall that \textbf{[Alg-2]} terminates at the stopping condition \textbf{[tiny-dom]}, i.e.
\begin{equation}\label{fs 1}
	\sum_{S_{m} \in \mathcal{S}_{m}}  \|e^{it\Delta}f_{S_{m} }\|^{p_{m}}_{BL^{p_{m},q}(B_{r_{m}})}  \leq 2\sum_{S_{m} \in \mathcal{S}_{m,\mathrm{tiny}}}  \|e^{it\Delta}f_{S_{m} }\|^{p_{m}}_{BL^{p_{m},q}(B_{r_{m}})}.
\end{equation}
For each $S_{m} \in \mathcal{S}_{m,\mathrm{tiny}} $, let $\mathcal{O}[S_{m}]$ denote the final collection of cells output by \textbf{[Alg-1]}. By {\bf Properties I, II} and {\bf III} of \textbf{[Alg-1]}, we obtain
\begin{equation}\label{fs 2}
	\|e^{it\Delta}f_{S_{m}}\|^{p_{m}}_{BL^{p_{m},q}(B_{r_{m}})} \lessapprox  D^{\delta}_{m-1} \sum_{O \in \mathcal{O}[S_{m}]}    \|e^{it\Delta}f_{O}\|^{p_{m}}_{BL^{p_{m},q}(O)},
\end{equation}
where the functions $f_{O}$ satisfy
\begin{equation}\label{fs 3}
	\sum_{O \in \mathcal{O}[S_{m}]}\|f_O\|_{L^2}^2 \lessapprox D_{m-1}^{1+\delta}\|f_{S_m}\|_{L^2}^2
\end{equation}
and
\begin{equation}\label{fs 4}
	\max_{O \in \mathcal{O}[S_{m}]}\|f_O\|_{L^2}^2 \lessapprox \Big( \frac{r_m}{r_{m-1}} \Big)^{-\frac{n+1-m}{2}}D_{m-1}^{-(m-1)+\delta}\|f_{S_m}\|_{L^2}^2.
\end{equation}
Let $\mathcal{O}$ denotes the union of the $\mathcal{O}[S_{m}]$ over all $S_{m}\in \mathcal{S}_{m,\mathrm{tiny}}$. Combining {\bf Property 1} of \textbf{[Alg-2]}, (\ref{fs 1}) and (\ref{fs 2}), one concludes
\begin{equation}\label{fs 5}
	\|e^{it\Delta}f\|_{BL^{p,q}(B^{\ast}_{R})} \lessapprox D^{\delta}_{m-1}M(\vec{r}_{m},\vec{D}_{m}) \|f\|^{1-\beta_m}_{L^{2}}\Big( \sum_{O \in \mathcal{O}}  \|e^{it\Delta}f_{O}\|^{p_{m}}_{BL^{p_{m},q}(O)}              \Big)^{\frac{\beta_{m}}{p_{m}}}  .
\end{equation}
Since each $O \in\mathcal{O}$ has diameter at most $ R^{\delta_{0}}$ by the stopping condition \textbf{[tiny]}, one has
$$ \|e^{it\Delta}f_{O}\|_{BL^{p_{m},q}(O)} \lessapprox \|f_{O}\|_{L^{2}}.$$
Hence \eqref{fs 5} becomes
\begin{equation}
	\|e^{it\Delta}f\|_{BL^{p,q}(B^{\ast}_{R})} \lessapprox D^{\delta}_{m-1}M(\vec{r}_{m},\vec{D}_{m}) \|f\|^{1-\beta_m}_{L^{2}}\Big( \sum_{O \in \mathcal{O}}  \|f_{O}\|^{p_{m}}_{L^2}              \Big)^{\frac{\beta_{m}}{p_{m}}}  .
\end{equation}
The definition of $\beta_{m}$ leads to
$$ \Big( \frac{1}{2}-\frac{1}{p_{m}}  \Big)\beta_{m}   = \frac{1}{2}-\frac{1}{p_{n+1}},$$
and so we obtain
\begin{equation}\label{fs 7}
	\|e^{it\Delta}f\|_{BL^{p,q}(B^{\ast}_{R})} \lessapprox D^{\delta}_{m-1}M(\vec{r}_{m},\vec{D}_{m}) \|f\|^{1-\beta_m}_{L^{2}}\Big( \sum_{O \in \mathcal{O}}  \|f_{O}\|^{2}_{L^2} \Big)^{\frac{\beta_{m}}{p_{m}}} \max_{O\in \mathcal{O}}\|f_O\|_{L^2}^{1-\frac{2}{p_{n+1}}}  .
\end{equation}

On the other hand, using (\ref{fs 3}) and {\bf  Property 2} from \textbf{[Alg-2]}, one has
\begin{equation}\label{fs 8}
	\sum_{O \in \mathcal{O}}\|f_O\|_{L^2}^2 \lessapprox  \Big(\prod_{i=m-1}^{n}D_i^{1+\delta} \Big) \|f\|_{L^2}^2  .
\end{equation}
By (\ref{fs 4}) and repeated application of {\bf Property 3} from \textbf{[Alg-2]}, it implies that
\begin{equation}\label{fs 9}
	\max_{O \in \mathcal{O}} \|f_{O}\|^{2}_{L^{2}} \lessapprox \prod_{i=m-1}^{n} r_{i}^{-\frac{1}{2}}D_{i}^{-i+\delta}\|f\|^{2}_{L^{2}},
\end{equation}
where $r_{m-1}:=1$.
Combining (\ref{fs 7}), (\ref{fs 8}) and (\ref{fs 9}), one concludes
\begin{equation*}
	\|e^{it\Delta}f\|_{BL^{p,q}(B^{\ast}_{R})} \lessapprox \prod_{i=m-1}^{n} r_{i}^{X_{i}}D_{i}^{Y_{i}+O(\delta)}\|f\|_{L^{2}},
\end{equation*}
where
\begin{equation}\label{Xi,Yi} \left\{\begin{aligned} &X_{i}:=\tfrac{i}{2(i+1)}(\beta_{i+1}-\beta_{i})-\tfrac{1}{2}\big(\tfrac{1}{2}-\tfrac{1}{p_{n+1}}\big),\\
		&Y_{i}:=\tfrac{\beta_{i+1}}{2}-(i+1)\big(\tfrac{1}{2}-\tfrac{1}{p_{n+1}}\big).  \end{aligned}\right.    \end{equation}
In order to prove  Theorem \ref{th4}, we need to verify
$$X_{i} \leq 0, \;\;\; m\leq i \leq n\; \;\;\;\text{and}\;\;\; \;Y_{i} \leq 0, \;\;\;  m-1 \leq i \leq n.$$
One easily verifies that
\begin{equation}\label{Xi define}
	X_{i}\leq 0 ~\Longleftrightarrow~ \Big( \frac{1}{2}-\frac{1}{p_{i+1}} \Big)^{-1}-\Big( \frac{1}{2}-\frac{1}{p_{i}} \Big)^{-1}\leq \frac{i+1}{i}
\end{equation}
and
\begin{equation}\label{Yi define}
	Y_{i-1}\leq 0 ~\Longleftrightarrow~ \Big( \frac{1}{2}-\frac{1}{p_{i}} \Big)^{-1}-2i\leq 0.
\end{equation}

We firstly choose $p_{m}=\frac{2m}{m-1}$, then
$$\Big( \frac{1}{2}-\frac{1}{p_{m}} \Big)^{-1} =2m.$$
Let $i=m$ in (\ref{Xi define}) and $i=m+1$ in (\ref{Yi define}), then
$$   \Big( \frac{1}{2}-\frac{1}{p_{m+1}} \Big)^{-1} \leq \min\Big\{ 2(m+1),2m+\frac{m+1}{m}  \Big\} =2m+\frac{m+1}{m}. $$
We choose $p_{m+1}$ such that
$$   \Big( \frac{1}{2}-\frac{1}{p_{m+1}} \Big)^{-1}  =2m+\frac{m+1}{m}. $$
Repeating the above argument, we obtain the following formula
$$   \Big( \frac{1}{2}-\frac{1}{p_{i}} \Big)^{-1}  =2m+\frac{m+1}{m}+...+\frac{i}{i-1},  \;\;\; i\geq m+1. $$
In particular, let  $i=n+1$,
$$ \Big( \frac{1}{2}-\frac{1}{p_{n+1}} \Big)^{-1}  =2m+\frac{m+1}{m}+...+\frac{n+1}{n},   $$
and so we have
$$  p_{n+1} =2+\frac{4}{n+m-1+\frac{1}{m}+...+\frac{1}{n}}.  $$
The worst case happens when $m=2$, i.e.
$$  p_{n+1} =2+\frac{4}{n+1+\frac{1}{2}+...+\frac{1}{n}}.  $$
Therefore, we have proved Theorem \ref{th4}.

\subsection*{Acknowledgements}    We  thank the anonymous referee and the associated editor  for their
invaluable comments   which helped to improve the paper. We are grateful for professor Xiaochun Li  warm comments and  suggestions. We are also grateful to Xiumin Du and Jianhui Li for pointing out an error in our initial proof of Theorem \ref{th4}.
This work is supported by the National Key research and development program of China [grant number 2020YFA0712903]; NSFC [grant number 11771388]; and NSFC [grant number 11831004].

\bibliographystyle{plain}
\bibliography{Cao-Miao-Wang}

\end{document}